\documentclass[12pt,a4paper]{book}

\usepackage{fancyhdr}
\usepackage{calc}

\fancypagestyle{plain}{
\fancyhead{} 
}
\makeatletter
\def\cleardoublepage{\clearpage\if@twoside \ifodd\c@page\else
\hbox{}
\thispagestyle{empty}
\newpage
\if@twocolumn\hbox{}\newpage\fi\fi\fi}
\makeatother

\usepackage[T1]{fontenc}

\usepackage{graphicx,psfrag,color,pst-grad}

\usepackage{amsfonts}
\usepackage{amsmath}
\usepackage{amssymb}
\usepackage{amsthm}
\usepackage{fancyvrb}
\usepackage{relsize}
\usepackage{graphicx}
\usepackage{textcomp}

\def\l2{\ell_2(X)}

\def\Rn{\mathbb{R}^n}
\def\R{\mathbb{R}}
\def\ns{{\mathcal{N}_{\Phi}(\Omega)} }
\def\hs{{\mathcal{N}_{\Phi}(X)}}
\def\T{{\mathcal{T}}}
\def\U{{\mathcal{U}}}
\def\L{{\mathcal{L}}}

\def\Uc{$\U  = \{u_i \in\ \ns,\ i = 1, \dots, N \} $}
\def\Vu{{V_{\U}}}
\def\Cu{{C_{\U}}}
\def\Vl{{V_{\L}}}
\def\Cl{{C_{\L}}}
\def\l{\ell_2(X)}
\def\lw{\ell_2^w(X)}
\def\C{{\mathcal{C}}}
\def\Pf{{\mathcal{P}_{\Phi,X}(x)}}
\def\Intx{{P_X[f](x)}}
\def\Int{{P_X[f]}}
\def\Appr{{\Lambda_{M}[f]}}
\def\Apprx{{\Lambda_{M}[f](x)}}
\def\T{{T_{\Phi}}}

\def\cubr{(X,\ \mathcal{W})_N}
\def\Mat{\textsl{Matlab}}
\def\cpp{\textsl{C$++$}}

\DeclareMathOperator*{\argmin}{arg\,min}

\newtheorem{theorem}{Theorem}[chapter]
\newtheorem{prop}[theorem]{Proposition}
\newtheorem{cor}[theorem]{Corollary}
\newtheorem{defin}{Definition}[chapter]
\newtheorem{remark}{Remark}[chapter]

\newenvironment{pf}
{\noindent \textbf{Proof:}}
{\qedsymbol\newline}

\begin{document}
\tableofcontents
\chapter*{Introduction}
\addcontentsline{toc}{chapter}{\numberline{}Introduction}
The main aim of Approximation Theory is to reconstruct a given function defined on a set $\Omega\subset\Rn$ from  
some values sampled on a finite set $X\subset\Omega$.
\\This process is required to be convergent and stable, namely to be such that under suitable conditions the approximant
is able to reproduce the original function with respect to a given norm.
\\\\In this setting, the so-called \textit{Kernel methods} are of growing importance. Altough they are built to 
be well-posed for every data distribution, it is also well-known that in many cases they suffers from serious instability
if no attention is paid to some aspects of their use.      
\\Several approaches have been studied to assure a fast convergence together with a stable computation. They are based
on different aspects of the approximation process, from the optimization of certain parameters and the research of convenient
data sets, to the numerical stabilization of the underlying linear system. 
\\Recently the two papers \cite{MuSch,SchPa} introduced a new tool, that is a general way to produce a stable basis for the 
functional space $\ns$ where the approximation takes place, based only on a particular factorization of the collocation 
matrix associated to the kernel. In the original works this process was used to create a \textit{Newton basis}, which is in
particular stable, complete, orthonormal and recursively computable. 
\\\\We use these results to build a different kind of stable basis that shares some useful properties with the one proposed
by the authors. Moreover, our basis provides a connection with a ``natural'' basis for the functional space $\ns$, which
arises from an eigendecomposition of a compact integral operator associated with the kernel, and which brings intrinsic 
information about the kernel itself and about the set $\Omega$. 
\\On the numerical side, the structure of the basis allows to further stabilize the approximation by moving from an exact 
data interpolation to an approximation in the least-squares sense, with a process that exactly 
corresponds to a low-rank approximation of the kernel matrix.    
\\\\The work is structured as follows: Chapter \ref{ch:general} gives a short introduction to the theory of 
Radial Basis Function; Chapter \ref{ch:change} describes the particular procedure introduced in the papers 
\cite{MuSch,SchPa} that will be the starting point for the development of our basis, which is constructed and analyzed in Chapter
\ref{ch:wsvd}; Chapter \ref{ch:numerics} presents some numerical examples which test our method from different 
point of view. The last Chapter discusses potential work that could be done to improve and better understand our results.

\chapter{Radial Basis Functions}\label{ch:general}
In this chapter we will give a brief introduction to the general theory of Radial Basis Functions (RBF).
\\We will start from some remarks on multivariate approximation that motivates the introduction of this kind of 
technique, then we will recall the basic theoretical results that arises in this context, focusing mainly on the tools that
will be used in the following chapters. 
\\The last section summarizes the main convergence and stability estimates, with particular focus on possible weakness of
RBF approximation.  
\\\\Each result presented here is taken from the book \cite{Bible}, that contains a complete coverage
of the theory of RBF as well as a much more general treatment of the Scattered Data Approximation. 

\section{Motivation}
The goal is to reconstruct a function $f\in\mathcal{C}(\Omega)$ defined on a set $\Omega\in\Rn$ from its samples at a 
fixed, discrete data-sites set $X\subset\Omega$, of cardinality $|X|=N\in\mathbb{N}$. 
\\The idea is to fix a finite dimensional subset $V\subseteq\mathcal{C}(\Omega)$ that allows a sufficiently good approximation
of the full space and to use the data-values $f_{|_{X}}$ to represent the function in this subspace.
\\To put this in practice, if V is spanned by a basis $\{v_1,\dots,v_N\}$, one wants to find an interpolant $P[f]\in V$ such that 
\begin{equation*}
P[f](x)=\sum_{j=1}^N c_j v_j(x)\quad \forall x\in\Omega,\quad P[f]_{|_{X}}=f_{|_{X}} 
\end{equation*}
where the coefficient vector $[c_1,\dots,c_N]^T$ is clearly the solution of the linear system
\begin{equation*}
\left( v_j(x_i) \right)_{i,j} \cdot c_j = f(x_i) \quad i,j=1,\dots,N
\end{equation*}
The matter is now to choose a good subset $V$, i.e. a subset for which the problem is well-posed. These space are the so-called 
\textit{Haar spaces}:
\begin{defin}[Haar space]
Suppose that $\Omega\subset\Rn$ contains at least N points. Let $V\subseteq \mathcal{C}(\Omega)$ be an $N$-dimensional linear space 
of functions, and let $\{v_1,\dots,v_N\}$ be a basis for $V$. Then $V$ is called a Haar space of dimension $N$ on $\Omega$ if
\begin{equation*}
\det(v_j(x_i))\neq0 
\end{equation*}
for any distinct points $x_1,\dots,x_N \in\Omega$. 
\end{defin}
\noindent This is, for example, the case of the space $\Pi_d$ of the univariate polynomial of degree at most $d$, 
that is indeed widely used and studied in approximation theory. 
\\Unfortunately if the dimension $n$ of the underlying real space is bigger than 1, there is no hope to find an Haar 
space:
\begin{theorem}[Haar - Mairhuber - Curtis]
Suppose that $\Omega\subseteq\Rn$, $n\geqslant2$, contains an interior point. Then there exists no Haar spaces on $\Omega$ of 
dimension $N\geqslant2$.
\end{theorem}
\noindent The last Theorem forces to use, in the multivariate setting, spaces and bases that depend on the choosen points.
\\\\In this context takes place the RBF approximation. Indeed, given a kernel 
$$\Phi:\Omega\times\Omega\rightarrow\R$$
we build a data-sites dependent basis  
\begin{equation}
\mathcal{T}_X = \{\Phi(\cdot,x_1),\dots,\Phi(\cdot,x_N) \} 
\end{equation}
that can be used to approximate a function $f$ as described.
\\In particular we will solve the linear system
\begin{equation*}
\left( \Phi_j(x_i,x_j) \right)_{i,j} \cdot c_j = f(x_i). 
\end{equation*}
The problem is well-posed if the so-called kernel matrix $A:=\left( \Phi_j(x_i,x_j) \right)_{i,j}$ is invertible for
every possible choice of the data-sites $X$.
\\The way to ensure this is to require that the kernel $\Phi$ is strictly positive definite: as shown by the following
definition, a positive definite kernel leads to a positive definite and hence invertible kernel matrix:
\begin{defin}
A continuous function $\Phi:\Rn\times\Rn\rightarrow\R$ is called positive definite if, for all $N\in\mathbb{N}$, for all sets of pairwise
distinct centers $X=\{x_1,\dots,x_N\}\in\Rn$ and for all $\alpha\in R^N\setminus\{0\}$, the quadratic form 
\begin{equation*}
\sum_{j=1}^N\sum_{i=1}^N\alpha_j\alpha_i\Phi(x_j,x_i) 
\end{equation*}
is positive.
\end{defin}
\noindent Even though this hypotesis on $\Phi$ suffices to have an unique solution for the approximation problem, a further requirement
gives some advantages: it's common to assume that $\Phi$ is radial, in the sense specified by the next definition.  
\begin{defin}
A function $\Phi:\Rn\times\Rn\rightarrow\R$ is said to be radial if there exist a function $\phi: [0,\infty)\rightarrow\R$ 
such that $\Phi(x,y)=\phi(\|x-y\|_2)$.
\end{defin}
\noindent This property allows to use the same structure for every dimension $n\in\mathbb{N}$ (dimension-blindness), and moreover
assures that the kernel is symmetric, together with the kernel matrix $A$. Nevertheless, there are some kernels that are
positive definite only for some dimension $n$ of the space.
\\\\In \cite[Ch.6]{Bible} is also possible to find a complete characterization of this kind of functions, based 
on generalized Fourier transforms.
\\Some relevant examples of radial and positive definite kernel functions that will be used in this work are listed in the table 
\eqref{tb:kernels}. It's common to use a radially scaled version 
\begin{equation*}
\Phi(x,y)=\phi(\varepsilon\|x-y\|_2) 
\end{equation*}
of the kernel, in order to have more control on the behaviour of the approximant. The parameter $\varepsilon\in\R$ is referred
as \emph{shape parameter}. 
\begin{table}
\centering
\begin{tabular}{||l|c|c||}
\hline 
&$\phi(r)$&dim\\
\hline
Gaussian&$e^{-(\varepsilon r)^2}$&$n$\\
&&\\
\hline
Inverse Multiquadric (IMQ)&$1/\sqrt{1+(\varepsilon r)^2} $&$n$\\
&&\\
\hline
Generalized IMQ&$1/(1+(\varepsilon r)^2)^2 $&$n$\\
&&\\
\hline
Inverse quadratic (IQ)&$1/(1+\varepsilon r) $&$n$\\
&&\\
\hline
Linear Mat\'ern (1MAT)&$e^{-\varepsilon r}(1+\varepsilon r)$&$n$\\
&&\\
\hline
Quadratic Mat\'ern (2MAT)&$e^{-\varepsilon r}(3+3 \varepsilon r + (\varepsilon r)^2)$&$n$\\
&&\\
\hline
Cubic Mat\'ern (3MAT)&$e^{-\varepsilon r}(15+15 \varepsilon r + 6 (\varepsilon r)^2+(\varepsilon r)^3))$&$n$\\
&&\\
\hline
Linear Laguerre-Gaussian&$e^{-(\varepsilon r)^2}(2-(\varepsilon r)^2))$&$2$\\
&&\\
\hline
Quadratic Laguerre-Gaussian&$e^{-(\varepsilon r)^2}(3-3(\varepsilon r)^2)+\frac{1}{2}(\varepsilon r)^4)$&$2$\\
&&\\
\hline
Linear generalized IMQ&$(2-(\varepsilon r)^2)/(1+(\varepsilon r)^2))^4$&$2$\\
&&\\
\hline
Wendland $2,0$ (W20)&$(1-\varepsilon r)_+^2 $&$2$\\
&&\\
\hline
Wendland $2,1$ (W21)&$(1-\varepsilon r)_+^4 (3 \varepsilon r + 1 )/20$&$2$\\
&&\\
\hline
\end{tabular} 
\label{tb:kernels}
\caption{Examples of radial positive definite kernels. The table shows the functions $\phi:\R\rightarrow\R$ such that 
$\Phi(\cdot,\cdot):= \phi(\|\cdot-\cdot\|_2)$ is positive definite in $\Rn$, where $n$ is as in the third column.}
\end{table}

\section{Native Space}\label{NativeSpace}
There is a natural space in which consider the RBF approximation. In fact, for each positive definite and symmetric kernel $\Phi$ 
and for each region $\Omega\subset\Rn$ it is possible to define an associated real Hilbert space, the so-called 
\emph{Native Space} $\ns$. In this space $\Phi$ is a reproducing kernel, in the sense of the following definition:
\begin{defin}[RK]
Let $\mathcal{H}$ be a real Hilbert space of functions $f: \Omega\rightarrow\R$. A function $\Phi: \Omega\times\Omega\rightarrow\R$
 is called a reproducing kernel for $\mathcal{H}$ if
\begin{enumerate}
\item $\Phi(\cdot,y)\in \mathcal{H}\quad \forall y\in\Omega$  
\item $f(y) = (f,\Phi(\cdot,y))_{\mathcal{H}}\quad \forall f\in\mathcal{H},\ \forall x\in\Omega$ (reproducing property)
\end{enumerate}
\end{defin}
\noindent It's well known that the existence of a positive definite reproducing kernel for an Hilbert space $H$ is equivalent 
to the continuity and the linear independence of the point evaluation functionals $\delta_x \in H^{*}$.  
\\\\The construction of such a space works as follows. Define the linear space
\begin{equation}\label{hs}
N_{\Phi}(\Omega) := \text{span}\{\Phi(\cdot,y):y\in\Omega\} 
\end{equation}
and equip it with the bilinear form
\begin{equation*}
\left(\sum_{j=1}^N \alpha_j \Phi(\cdot,x_j),\sum_{i=1}^M \beta_i \Phi(\cdot,x_i)\right)_{\Phi}:=\sum_{j=1}^N\sum_{i=1}^M \alpha_j \beta_i \Phi(x_j,x_i) .
\end{equation*}
Then the space $\left(N_{\Phi}(\Omega), (\cdot,\cdot)_{\Phi} \right)$ is almost what we want:
\begin{theorem}
Under the above assumptions on the kernel, $(\cdot,\cdot)_{\Phi}$ defines an inner product on 
$N_{\Phi}(\Omega)$, which is a pre-Hilbert space with reproducing kernel $\Phi$.
\end{theorem}
\noindent The completion of $ N_{\Phi}(\Omega)$ with respect to the $\|\cdot\|_{\Phi}$-norm is the Native Space $\ns$, which 
is an Hilbert space with $\Phi$ as reproducing kernel. 
\\This space clearly contains elements $f$ not in $N_{\Phi}(\Omega)$, that can 
be understood as functions of the form
\begin{equation*}
f(x) = (f,\Phi(\cdot,x))_{\Phi} \quad\forall x\in\Omega
\end{equation*}
since the continuity of $\delta_x$ is preserved in the completion. 
\\\\This space is unique in the sense that if $\mathcal{G}$ is a real Hilbert space of functions $f:\Omega\rightarrow\R$
with reproducing kernel $\Phi$, then $\mathcal{G}=\ns$ and the inner products are the same. 
\\Moreover it can be shown that if an Hilbert space $\mathcal{G}$ has a reproducing kernel, it's necessarily positive 
definite and symmetric, and then there is a full equivalence between Reproducing Kernel Hilbert Spaces and Native Spaces for
positive definite and symmetric kernels.
\subsection{Embeddings}
For their use in approximation, it is useful to know what kind of functions belongs to the Native Space and how the smoothness 
of the kernel is inherited: 
\begin{theorem}[$\mathcal{C}^k(\Omega)$]\label{embed1}
Let $\Omega\subset\Rn$ be an open set and let $\Phi\in\mathcal{C}^{2k}(\Omega\times\Omega)$ be a symmetric, positive 
definite and radial kernel on $\Omega$. Then $\ns \subset \mathcal{C}^{k}(\Omega)$ and $\forall \alpha\in\mathbb{N}^n$, 
$|\alpha|\leqslant k$, $\forall f\in\ns$ and $\forall x\in\Omega$ 
\begin{equation*}
D^{\alpha}f(x) = (f,D_2^{\alpha}\Phi(\cdot,x))_{\Phi}
\end{equation*}
\end{theorem}
\noindent The case k=0 in particular ensures that the $\|\cdot\|_{\Phi}$-convergence implies the pointwise convergence
(it's also a direct consequence of the reproducng kernel property). 
\begin{theorem}[$L_2(\Omega)$]\label{th:embed2}
Let $\Omega\subset\Rn$ be a compact set and let $\Phi$ be a symmetric, positive definite and radial kernel on $\Omega$. Then
the Native Space $\ns$ has a continuous linear embedding into $L_2(\Omega)$, and in particular
\begin{equation*}
\|f\|_{L_2(\Omega)}\leqslant \sqrt{|\Omega|\ \phi(0) }\  \|f\|_{\Phi}\quad \forall f \in \ns
\end{equation*}
where $|\Omega|:=\text{meas}(\Omega)$ 
\end{theorem}
 
\subsection{An integral operator and a ``natural'' basis}
It is possible to define an integral operator $T_{\Phi}$ associated to the kernel and defined on the Native Space.
\\Construction and properties of such operator are discussed in detail in the paragraph 10.4 of the book \cite{Bible}, where it is the base 
for a further characterization of $\ns$. Here we are interested mainly in a particular basis that arises from an 
eigendecomposition of it, which would be the key tool for the subsequent discussion.
\\Consider the operator $T_{\Phi} : L_2(\Omega) \rightarrow L_2(\Omega)$ defined by
\begin{equation}\label{theOperator}
T_{\Phi}[f] (x) := \int_{\Omega} \Phi(x,y) f(y) dy\quad \forall f\in L_2(\Omega),\ \forall x\in\Omega 
\end{equation}
that maps $L_2(\Omega)$ continuously into $\ns$. It is the adjoint of the embedding operator of 
$\ns$ into $L_2(\Omega)$, i.e.
\begin{equation}\label{L2N}
(f,v)_{L_2(\Omega)} = (f,T_{\Phi}[v])_{\Phi} \quad \forall f\in\ns,\ \forall v\in L_2(\Omega).
\end{equation}
A particular and in some sense ``natural`` basis for $\ns$ comes from the following theorem:
\begin{theorem}[Mercer]\label{th:mercer}
Every continuous positive definite kernel $\Phi$ on a bounded domain $\Omega\subset\Rn$ defines an operator
\begin{equation*}
\T: \ns \rightarrow \ns, \quad \T[f] = \int_{\Omega} \Phi(x,y) f(y) dy 
\end{equation*}
which is bounded, compact and self-adjoint. It has an enumerable set of eigenvalues and eigenvectors 
$\{\varphi_j\}_{j>0}$,
\begin{equation*}
\lambda_j \varphi_j(x) = \int_{\Omega} \Phi(x,y) \varphi_j(y) dy\quad \forall x\in\Omega\label{def:operator}\\
\end{equation*}
which forms an orthonormal basis for $\ns$, and in particular
\begin{eqnarray*}
\{\varphi_j\}_{j>0} &&is\ orthonormal\ in\ \ns \\
\{\varphi_j\}_{j>0} &&is\ orthogonal\ in\ L_2(\Omega),\ \|\varphi_j\|_{L_2(\Omega)}^2= \lambda_j\\
\lambda_j\rightarrow0 &&as\ j\rightarrow\infty
\end{eqnarray*}
Moreover the kernel has an eigenfunctions expansion
\begin{equation*}
\Phi(x,y)=\sum_{j=1}^{\infty}\lambda_j\ \varphi_j(x) \varphi_j(y)\quad  \forall x,y \in \Omega 
\end{equation*}
which is absolutely and uniformly convergent. 
\end{theorem}
\begin{remark}\label{rem:trace}
The operator $T_{\Phi}$ is a trace-class operator, and in particular
\begin{equation*}
\sum_{j>0} \lambda_j = \int_{\Omega} \Phi(x,x)\ dx = \phi(0)\ |\Omega| 
\end{equation*}
This property, together with the fact that the eigenvalues accumulates in $0$, will be useful to estimate
the convergence of the truncated series to the full one.
\end{remark}
\begin{remark}\label{L2Nphi}
A consequence of the property \eqref{L2N}, which we point out for later use, is that $\forall j>0$
\begin{equation*}
(f,\varphi_j)_{L_2(\Omega)} = (f,T_{\Phi}[\varphi_j])_{\Phi} = \lambda_j\ (f,\varphi_j)_{\Phi} = 
(\varphi_j,\varphi_j)_{L_2(\Omega)}\ (f,\varphi_j)_{\Phi} \quad \forall f\in\ns
\end{equation*}
\end{remark}

\subsection{Other inner products}
For later use we introduce also the following discrete scalar products, which will be of key importance in the following 
chapters.
\begin{defin}[$\l$]\label{l2prod}
Let $f$, $g$ be functions in $\ns$ and $X=\{x_1,\dots,x_N\}\subset\Omega$ a discrete set. The $\l$-scalar product is defined as 
\begin{equation*}
(f,g)_{\l} = \sum_{i=1}^N f(x_i) g(x_i) 
\end{equation*}
\end{defin}
 
\begin{defin}[$\lw$]\label{l2wprod}
Let $f$, $g$ be functions in $\ns$, $X=\{x_1,\dots,x_N\}\subset\Omega$ a discrete set and $\mathcal{W}=\{w_1,\dots,w_N\}\subset\R$
a set of positive weights. The $\lw$-scalar product is defined as 
\begin{equation*}
(f,g)_{\lw} = \sum_{i=1}^N w_i\ f(x_i) g(x_i) 
\end{equation*}
\end{defin}
\noindent It is clear from the definitions that both the products are not positive definite, since for each $f\in\ns$ that 
vanishes on $X$, $(f,f)_{\l}=(f,f)_{\lw}=0$. Anyway they are positive definite when restricted to $\hs$.

\section{Error bounds and stability estimates}\label{sec:stdconvstab}
We recall that, for a subset $\Omega\subset\Rn$, a discrete data-sites set $X\subset\Omega$ and a radial, positive 
definite kernel $\Phi\in\mathcal{C}(\Omega\times\Omega)$ the RBF interpolant $\Int$ to a function $f\in\ns$ is computed as
\begin{equation}\label{interpolant0}
\Intx=\sum_{j=1}^N c_j\ \Phi(x,x_j),\quad \Int(x_i)=f(x_i)\quad \forall x\in\Omega, x_i\in X
\end{equation}
The question is how well $\Int$ can approximate the sampled function $f$, i.e. if $\Int$ converges to $f$ in some 
given norm when the data-sites $X$ becomes dense in $\Omega$. To be more precise, one wants to know
if $\hs$ saturates $\ns$ for a good choice of $X\subset\Omega$, and if the process can be accomplished in a stable way.
\\\\There are two quantities used to relate the set $X$ to these requirements: the \textit{fill distance}
\begin{equation}
h_{X,\Omega} = \max_{x\in\Omega} \min_{x_i\in X} \|x-x_i\|_2
\end{equation}
and the \textit{separation distance}
\begin{equation}
q_{X} = \frac{1}{2}\min_{x_i,x_j\in X} \|x_j-x_i\|_2
\end{equation}
Clearly the shape parameter $\varepsilon$ have also an important role, since it determines the radial amplitude 
of the kernel.
\\\\The first estimate comes directly from the definition of the pointwise-error functional. Let $\mathcal{E}_x$
be defined $\forall x\in\Omega$ as
\begin{equation*}
\mathcal{E}_x: \ns\rightarrow \R,\quad \mathcal{E}_x[f] = f(x)-\Intx
\end{equation*}
and let $\mathcal{P}_{\Phi,X}$ denote its norm, the so-called \textit{Power Function}.
\\Then the basic estimate for the convergence is the following:
\begin{theorem}\label{th:bound1}
Let $\Omega\subset\Rn$, let $\Phi\in\mathcal{C}(\Omega\times\Omega)$ be a symmetric positive definite kernel, let $X\subset\Omega$
be a discrete set of data sites and let $f\in\ns$. Then
\begin{equation}
|f(x)-\Intx| \leqslant \Pf \|f\|_{\Phi}\quad \forall x\in\Omega
\end{equation}
\end{theorem}
\noindent It can be refined observing that the interpolation operator is a projection with respect to the 
$\Phi$-inner product (and then the interpolant is also the best approximation in $\hs$ of $f\in\ns$). The estimate
then becomes
\begin{equation}
|f(x)-\Intx| \leqslant \Pf \|f-\Int\|_{\Phi}\quad \forall x\in\Omega
\end{equation}
Moreover, the Power Function can be exactly computed introducing a Lagrange basis for $\hs$:
\begin{prop}\label{def:Lb}
For any pairwise distinct data-sites set $X\in\Omega$ there exist a Lagrange basis $\L=\{\ell_1,\dots,\ell_N\}$ for $\hs$.
This basis is such that $\ell_j(x_i)=\delta_{ij}\ \forall i,j=1,\dots,N $ and then the interpolant can be written in 
cardinal form as
\begin{equation*}
\Intx=\sum_{j=1}^N f(x_j) \ell_j(x)\quad\forall x\in\Omega 
\end{equation*}
\end{prop}
\noindent The expression for the interpolant gives
\begin{equation}\label{def:Pf}
\Pf^2=\Phi(x,x)-\sum_{j=1}^N \Phi(x,x_j) \ell_i(x)
\end{equation}
Using this explicit representation it is possible to bound the Power Function on domains that satisfies an interior cone 
condition, using a multivariate Taylor expansion. This final estimate relates the set $X$ to the interpolation error: 
\begin{theorem}
Let $\Omega\subset\Rn$ be a bounded set that satisfies an interior cone condition and let 
$\Phi\in\mathcal{C}^{2k}(\Omega\times\Omega)$ be a symmetric positive definite kernel. Then there exist
positive constants $h_0$ and $C$ independent of $x$, $f$ and $\Phi$, such that $\forall X\subset\Omega$, 
$h_{X,\Omega}\leqslant h_0$, $\forall f\in\ns$ and $\forall x \in \Omega$ 
\begin{equation*}
|f(x)-\Intx| \leqslant C\cdot h_{X,\Omega}^k\cdot C_{\Phi} \|f-\Int\|_{\Phi}
\end{equation*}
where $C_{\Phi}$ is a constant that depends on the derivatives of $\Phi$.  
\end{theorem}
\noindent From the previous estimate we can expect that the approximation error goes to zero as $h_{X,\Omega}\rightarrow 0$. This 
is not completely true, since when the data-sites set $X$ becomes too big the interpolation can be instable. 
\\In fact, it is possible to prove that the condition number of the kernel matrix $A$ grows if the separation distance 
$q_{X}$ decreases, and this, together with a bad choice of the shape parameter $\varepsilon$, can produce very instable 
approximants. 
\\Various approaches are used to avoid this situation. A lot of efforts are made on the study of 
well-distributed data-sites set, for examples sets $X$ such that the \textit{uniformity}
\begin{equation*}
\rho_{X,\Omega} = \frac{q_{X}}{h_{X,\Omega}} 
\end{equation*}
is maximized.   
\\Another common way to try to avoid instability, and more related on the liner algebra part of the method, is to choose
a shape parameter $\varepsilon$ such that the kernel matrix is not ill-conditioned.
\\Recently another method to ensure convergence and stability was presented, and it is described in the next chapter.

\chapter{General bases}\label{ch:change}
As shown in the previous chapter, the use of the standard basis of translates leads to ill-conditioned kernel matrices, and, 
moreover, it gives poor information about the selection of a ``good'' centers set $X\subset\Omega$.  
\\Hence it makes sense to consider different bases of $\ns$ in order to obtain better results in term of stability and 
convergence, and such that it will be possible to describe in an useful way the data dependence of the subspace $\hs$.
\newline In the present chapter we will describe this kind of change of basis, focusing mainly on $\Phi$-orthonormal basis.
\\\\This approach was introduced in the papers \cite{MuSch, SchPa}, where it is the starting point to produce a 
\textit{Newton basis}.
\section{Definition and characterization} 
Let $\Omega \subset \Rn$, $X = \{x_1, \dots, x_N \}\subset\Omega$, and let $\mathcal{T}_X=\{\Phi(\cdot,x_i),x_i\in X\}$  be the standard basis of translates. 
\\Consider another basis {\Uc} such that 
\begin{equation}\label{basis}
\mathrm{span}(\U)= \mathrm{span}(\mathcal{T}_X ) =  \hs 
\end{equation}
The following Theorem gives a characterization of such bases, where $\mathcal{T}_X$ and $\U$ are expressed as row vectors: 
\begin{eqnarray*}
T(x) &=& [\Phi(x,x_1), \dots, \Phi(x,x_N) ] \in \R^N \\
U(x) &=& [u_1(x), \dots, u_N(x) ] \in \R^N 
\end{eqnarray*}

\begin{theorem}[Characterization]\label{th:char}
Any basis $\U$ arises from a factorization of the kernel matrix
\begin{equation}\label{A}
A=\Vu\cdot {\Cu}^{-1} 
\end{equation}
where 
\begin{equation}\label{Vu}
\Vu = (u_j(x_i))_{1\leqslant i,j\leqslant N}
\end{equation}
and the coefficient matrix $\Cu$ is such that 
\begin{equation}\label{Cu}
U(x) = T(x)\cdot \Cu
\end{equation}
\end{theorem}
\begin{pf}
Let $\U$ be a basis as defined in (\ref{basis}).
Since $u_i \in\ \hs\ \forall\ u_i \in \U$, there exist real coefficients $(c_{j i})_{1\leqslant i,j\leqslant N} $ 
such that $\forall x \in \Omega$
\begin{equation}\label{ui}
u_i(x) = \sum_{j=1}^N \Phi(x,x_j)\ c_{j i}\ , \:  1\leqslant i \leqslant N
\end{equation}
and it suffices to set $\Cu = (c_{j i})_{1\leqslant i,j\leqslant N} $ to obtain (\ref{Cu}).
\\Now consider the evaluation operator $E_X: \ns \rightarrow \R^N$,
\begin{equation*}
E_X(f) = [f(x_1), \dots, f(x_N) ]
\end{equation*}
that maps functions into column vectors and rows of functions into matrices. By definition
\begin{eqnarray*}
E( T(x)) &=& (\Phi(x_i,x_j))_{1\leqslant i,j \leqslant N} = A\\
E(U(x)) &=& (u_j(x_i))_{1\leqslant i,j \leqslant N}.
\end{eqnarray*}
Setting $\Vu := E(U(x))$ as in (\ref{Vu}), by (\ref{Cu}) we have 
\begin{equation*}
\Vu = E(U(x)) = E(T(x) \cdot \Cu)= A\cdot \Cu
\end{equation*}
and hence $A=\Vu\cdot {\Cu}^{-1}$. 
\end{pf}

\noindent Using this characterization is simple to describe useful properties of the basis $\U$. In particular we will be 
interested in bases that are orthogonal or even orthonormal with respect to the inner products defined in section 
\eqref{NativeSpace}. The gramians can be computed as follows: 

\begin{prop}[gramians]\label{th:gram}
Given a basis $\U$ as in (\ref{basis}), we have
\begin{eqnarray*}
G_{\U} &:=& ((u_i,u_j)_{\Phi})_{1 \leqslant i,j \leqslant N} = \Cu^T\cdot A\cdot\Cu \\
\Gamma_{\U} &:=& ((u_i,u_j)_{\l})_{1 \leqslant i,j \leqslant N} = \Cu^T\cdot A^2 \cdot\Cu\\
\Gamma_{\U}^w &:=& ((u_i,u_j)_{\lw})_{1 \leqslant i,j \leqslant N} = \Vu^T\cdot W\cdot \Vu 
\end{eqnarray*}
\end{prop}
\begin{pf}
The formula for the $\Phi$-gramian $G_{\U}$ comes from \eqref{ui}:
\begin{eqnarray*}
G_{\U} &=& ((u_i,u_j)_{\Phi})_{1 \leqslant i,j \leqslant N} 
       = \left( \sum_{h,k=1}^N c_{ki} c_{hj} \Phi(x_h,x_k)\right)_{1 \leqslant i,j \leqslant N} \\
       &=& \Cu^T\cdot A\cdot\Cu 
\end{eqnarray*}
The $\l$ and $\lw$-gramians $\Gamma_{\U}$ and $\Gamma^w_{\U}$ can be directly computed using \eqref{A} and \eqref{Vu}:
\begin{eqnarray*}
\Gamma_{\U} &:=& \left((u_i,u_j)_{\l}\right)_{1 \leqslant i,j \leqslant N} 
	     = \left( \sum_{h=1}^N u_i(x_h) u_j(x_h)\right)_{1 \leqslant i,j \leqslant N} \\
	    & =& \Vu^T\cdot \Vu = \Cu^T\cdot A^2 \cdot\Cu\\ 
\Gamma^w_{\U} &:=& \left((u_i,u_j)_{\lw}\right)_{1 \leqslant i,j \leqslant N} 
	     = \left( \sum_{h=1}^N w_h u_i(x_h) u_j(x_h)\right)_{1 \leqslant i,j \leqslant N} \\
	    & =& \Vu^T\cdot W \cdot \Vu  
\end{eqnarray*}
This concludes the proof.
\end{pf} 

\noindent For later use it is convenient to rewrite also the Lagrange basis (\ref{def:Lb}) and the Power Function (\ref{def:Pf}) using
the same notation:
\begin{prop}[Lagrange basis]
The Lagrange basis $\L$ for $\hs$ is described by the matrices 
\begin{equation*}
\Vl = I,\;\Cl=A^{-1},\;G_{\L}=A^{-1},\;\Gamma_{\L}=I
\end{equation*}
\end{prop}
\begin{pf}
The statement is a direct consequence of Theorem \eqref{th:char} and Proposition \eqref{th:gram} applied to the definition \eqref{def:Lb} of $\L$.
\end{pf} 

\begin{prop}[Power function]\label{char:pf}
The Power Function can be expressed as
\begin{eqnarray}\label{th:power}
\Pf &=& \Phi(x,x) - U(x)\cdot {G_{\U}}^{-1}\cdot U^T(x)\\
    &=& \phi(0) - U(x)\cdot {G_{\U}}^{-1}\cdot U^T(x) \quad \forall x \in \Omega
\end{eqnarray}
\end{prop}
\begin{pf}
The last Proposition gives 
\begin{equation*}
L(x) := [\ell_1(x),\dots,\ell_N(x)] = T(x)\cdot \Cl =  T(x) \cdot A^{-1}  
\end{equation*}
Now starting from the definition (\ref{def:Pf}) we get
\begin{eqnarray*}
\Pf &=& \Phi(x,x) -\sum_{i=1}^N \Phi(x,x_i) \ell_i(x) \\
&=& \Phi(x,x) - T(x)\cdot L^T(x)\\
&=& \Phi(x,x) - T(x)\cdot A^{-1} \cdot T^T(x) \\
&=& \Phi(x,x) - U(x)\cdot \Cu^{-1}\cdot A^{-1}\cdot (\Cu^{-1})^T \cdot U^T(x) \\
&=& \Phi(x,x) - U(x)\cdot {G_{\U}}^{-1}\cdot U^T(x)\\
&=& \phi(0) - U(x)\cdot {G_{\U}}^{-1}\cdot U^T(x)
\end{eqnarray*}
where we used the definitions of $U(x)$ and $G_{\U}$.
\end{pf} 

\section{Interpolation and stability}
Now it is possible to express the interpolant to a given function $f\in\ns$ using this notation.
\begin{prop}[Interpolation] \label{prop:int}
The interpolant $\Int$ to a function $f\in\ns$ on $X\subset\Omega$ can be rewritten as 
\begin{equation}\label{th:int}
\Intx = \sum_{j=1}^N \Lambda_j(f)\ u_j(x)= U(x) \cdot \Lambda_{\U}(f) \quad \forall x \in \Omega  
\end{equation}
where $\Lambda_{\U}(f) = [\Lambda_1(f), \dots, \Lambda_N(f)]^T\in\R^N$ is a column vector of values of linear functionals defined by
\begin{equation}
\Lambda_{\U}(f) = \Cu^{-1}\cdot A^{-1}\cdot E_X(f)=\Vu^{-1} E_X(f)
\end{equation}
\end{prop}
\begin{pf}\
The interpolant \eqref{interpolant0} is espressed by
\begin{equation}
\Intx = \sum_{j=1}^N\ \alpha_j \Phi(x,x_j) = T(x) \cdot \alpha 
\end{equation}
where $\alpha\in\R^N,\; A\alpha=E_X(f)$. Thus, according to \eqref{th:char},
\begin{eqnarray}
\Intx &=&  T(x) \cdot \alpha =  U(x)\cdot\Cu^{-1}\cdot A^{-1}\cdot E_X(f) \\
      &=& U(x)\cdot \Vu^{-1} E_X(f)
\end{eqnarray}
This proves the statement.
\end{pf}

\begin{prop}[Stability] 
Let $\rho(\cdot)$ be the spectral radius and let $\kappa_2(\cdot)$ be the condition number with respect to the 
euclidean norm ${\|\cdot\|}_2$ on $\R^N$. Then the stability of the evaluation of $\Int$ can be bounded 
$\forall x\in \Omega$ as
\begin{equation}
|\Intx|^2\leqslant {\|U(x)\|}_2^2\ {\|\Lambda_{\U}(f)\|}_2^2\leqslant\kappa_2(G_{\U})\  \phi(0)\ {\|f\|}_{\Phi}^2
\end{equation}
and in particular the following bounds hold:
\begin{eqnarray}
{\|U(x)\|}_2^2&\leqslant&\rho(G_{\U})\ \phi(0)\quad \forall x\in\Omega\\
{\|\Lambda_{\U}(f)\|}_2^2 &\leqslant& \rho(G_{\U}^{-1})\ {\|f\|}_{\Phi}^2\quad \forall f\in\ns
\end{eqnarray}
\end{prop}
\begin{pf}
By the previous Proposition and the H\"older inequality we get the bound
\begin{equation}
|\Intx|\leqslant {\|U(x)\|}_2\ {\|\Lambda_{\U}(f)\|}_2\quad\forall x\in \Omega
\end{equation}
The term ${\|U(x)\|}_2$ can be bounded using \eqref{th:power}: the Power Function is non negative and hence 
$U(x)\cdot {G_{\U}}^{-1}\cdot U^T(x) \leqslant \Phi(x,x)= \phi(0)$. Now, using the properties of the Rayleigh quotient 
for the symmetric matrix $G_{\U}^{-1}$, we get the desired bound. 
\\Now consider the term ${\|\Lambda_{\U}(f)\|}_2$: from \eqref{th:int} it follows that
\begin{eqnarray}\label{lambda}
{\|\Int\|}_{\Phi} &=& \sum_{i,j=1}^N \Lambda_j(f)\ \Lambda_i(f)\ (u_j(x),u_i(x))_{\Phi}\\
		  &=& \Lambda_{\U}^T\cdot G_{\U}\cdot\Lambda_{\U} 
\end{eqnarray}
and then the bound holds thanks to ${\|\Int\|}_{\Phi}\leqslant{\|f\|}_{\Phi}$, and applying the same eigenvalue 
manipulation as above.
\end{pf}

\begin{cor}[$\Phi$-orthonormal bases]
If $\U$ is a $\Phi$-orthonormal basis, the stability estimate becomes
\begin{eqnarray}
\left| \Intx \right| \leqslant \sqrt{\phi(0)}\ \|f\|_{\Phi} \quad \forall x\in \Omega
\end{eqnarray}
In particular, the value of $\|U(x)\|_2$ for fixed $x\in\Omega$ and the value of $\|\Lambda(f)\|_2$ for fixed $f\in\ns$
are the same for all $\Phi$-orthonormal basis, and the following bounds, that are not dependent on $X\subset\Omega$, 
hold $\forall x \in \Omega$ 
\begin{equation}
\|U(x)\|_2 \leqslant\sqrt{\phi(0)},\quad \|\Lambda(f)\|_2\leqslant \|f\|_{\Phi}
\end{equation}
\end{cor}
\begin{pf}
The bounds follow immediately from the previous Proposition in the case $G_{\U}=I$.  
\\The equality $G_{\U}=I$ and the Power Function equation \eqref{def:Pf} gives also,  $\forall x\in\Omega$,
\begin{equation}
\|U(x)\|_2 = U(x)\cdot U(x)^T = \Phi(x,x) - \Pf\leqslant \Phi(x,x)  
\end{equation}
and since $\Pf$ is independent of the particular basis chosen, so is the value ${\|U(x)\|}_2$. 
\\Finally \eqref{lambda} proves that $\|\Lambda(f)\|_2$ is independent on $\U$ if it is $\Phi$-orthonormal.    
\end{pf}

\section{Orthonormal bases}  
The previous Corollary suggests the opportunity of using orthonormal bases. 
The following Theorems give a complete characterization of bases $\U$ as defined in \eqref{basis} that are orthonormal with 
respect to the considered inner products. 
\begin{theorem}[$\Phi$-orthonormal bases]
Each $\Phi$-orthonormal basis $\U$ arises from a decomposition 
\begin{equation}
A=B^T\cdot B, 
\end{equation}
with $\Vu=B^T$ and $\Cu=B^{-1}$  
\end{theorem}
\begin{pf}
Since $G_{\U} =  \Cu^T\cdot A\cdot\Cu $ by Theorem \eqref{th:gram}, the condition $G_{\U}=I$ is equivalent to 
$A=(\Cu^T)^{-1} \cdot \Cu^{-1}$, and the statement holds according to th.\eqref{th:char}.    
\end{pf}
\\It is also useful for the next chapter to characterize bases that are orthonormal with respect to the discrete scalar product
introduced in chapter \ref{ch:general}:
\begin{theorem}[$\l$-orthonormal bases]
Each $\l$-orthonormal basis $\U$ arises from a decomposition 
\begin{equation}
A=Q\cdot B,\;  Q^T\cdot Q= I
\end{equation}
with $\Vu=Q$ and $\Cu=B$.  
\end{theorem}
\begin{pf}
Since $\Gamma_{\U} =  \Cu^T\cdot A^2\cdot\Cu $ by Theorem \eqref{th:gram}, the condition $\Gamma_{\U}=I$ implies that 
$A\cdot\Cu=\Vu$ is orthogonal.   
\end{pf}

\begin{remark}\label{rem:intorth}
If the basis $\U$ is $\Phi$-orthonormal, the functionals $\Lambda_j$ in \eqref{prop:int} are obviously $\Phi$-scalar 
products, and so the interpolant takes the form
\begin{equation}
\Int(x) = \sum_{j=1}^N (f,u_j)_{\Phi} u_j(x)
\end{equation}
indeed
\begin{eqnarray*}
\Lambda_j(f) &=&  (\Vu^{-1}\cdot E_X(f))_j = (\Cu^T\cdot E_X(f))_j = \sum_{i=1}^N c_{ij} f(x_i) \\
&=&\left(\sum_{i=1}^N c_{ij} \Phi(x_i,\cdot) ,f\right)_{\Phi} = (u_j,f)_{\Phi}
\end{eqnarray*}
The above fact is clear since the interpolation operator is a projection operator on $\hs$ with respect to 
$(\cdot,\cdot)_{\Phi}$, and $\U$ is a $\Phi$-orthonormal basis for $\hs$. 
\end{remark}

\chapter{Weighted SVD bases}\label{ch:wsvd}
This chapter deals with the main issue of the present work.
\\We will introduce a particular basis for the native space $\ns$, and we will present the reason to threat it.
In particular the discussion points out the stability and convergence rate of the interpolant and the 
weighted least-squares approximant based on it.
\\\\The main idea is to discretize the ``natural'' bases described in Theorem \ref{th:mercer}. 
\\The numerical approximation of 
such basis gives a point-dependent, discrete basis which can be described using the notations introduced in the previous 
chapter. 
\\The interest on this basis is that it preserves the properties of the continuous basis in a discrete setting, and 
in particual is a complete, $\Phi$-orthonormal and $\lw$-orthogonal basis for $\ns$.
\\\\The idea for using this approach comes again from the paper \cite{SchPa}, where it is suggested as a possible way to 
create a connection between the continuous basis and the discrete one. 
\section{Definition and basic properties}
Consider a cubature rule $ \cubr $, $N\in\mathbb{N}$, on $\Omega$, i.e. a set of points $X=\{x_j\}_{j=1}^N \subset \Omega$ and a set of positive
weights $\mathcal{W}=\{w_j\}_{j=1}^N$ such that
\begin{equation}
\int_{\Omega} f(y) dy \approx \sum_{j=1}^N f(x_j) w_j\quad \forall f\in\ns
\end{equation}
We can approximate the operator \eqref{def:operator} for each eigenvalue $\lambda_j$ using the Nystr\"om method based on 
the above cubature rule. A complete description of the method is present in \cite{Atk}.
\\The operator can be evaluated on $X$,
\begin{equation*}
\lambda_j \varphi_j(x_i) = \int_{\Omega} \Phi(x_i,y) \varphi_j(y) dy\quad\forall i=1,\dots,N,\ \forall j>0
\end{equation*}
and then discretized using the cubature rule as
\begin{equation}\label{eq:nystrom}
\lambda_j \varphi_j(x_i) = \sum_{h=1}^N \Phi(x_i,x_h) \varphi_j(x_h) w_h\quad\forall i,j=1,\dots,N
\end{equation}
Now, setting $W=diag(w_j)$, it suffices to solve the following discrete eigenvalue problem in order to find the desired 
approximation of the operator's eigenvalues and eigenfunctions (evaluated on $X$): 
\begin{equation*}
\lambda v = (A \cdot W) v.
\end{equation*}
However this approach does not lead directly to a connection between the discretized version of the basis of Theorem 
\ref{th:mercer} and a basis of the subspace $\hs$ as defined in equation \eqref{hs}, since it involves a scaled version $A\cdot W$ 
of the kernel matrix which is no more symmetric and that cannot be described as a factorization of $A$, as required by the construction made in the
previous chapter.
\\\\A solution is to rewrite \eqref{eq:nystrom} using the positivity of the weights as
\begin{equation}
\lambda_j (\sqrt{w_i} \varphi_j(x_i)) = \sum_{h=1}^N (\sqrt{w_i} \Phi(x_i,x_h) \sqrt{w_h}) (\sqrt{w_h} \varphi_j(x_h)) 
\quad \forall i,j=1,\dots,N
\end{equation}
and then to consider the corresponding eigenvalue problem
\begin{equation*}
\lambda \left(\sqrt{W}\cdot v\right) = \left(\sqrt{W}\cdot A \cdot \sqrt{W}\right) \left(\sqrt{W}\cdot v\right)
\end{equation*}
which is equivalent to the previous one, now involving the symmetric and positive definite matrix 
$A_W:=\sqrt{W}\cdot A \cdot \sqrt{W}$. 
\\In particular this matrix is normal, and then a singular value decomposition of it 
is also an unitary diagonalization.
\\\\Motivated by this approach we can introduce a new basis for $\hs$, described in terms of the notation given in
Theorem \ref{th:char}: 
\begin{defin}[Weighted Svd basis]
A weighted svd basis $\U$ is a basis for $\hs$ characterized by the following matrices:
\begin{equation*}
\Vu = \sqrt{W^{-1}}\cdot Q \cdot \Sigma ,\;\Cu=\sqrt{W}\cdot Q \cdot \Sigma^{-1} 
\end{equation*}
where 
\begin{equation*}
\sqrt{W}\cdot A \cdot \sqrt{W} = Q \cdot \Sigma^2 \cdot Q^T  
\end{equation*}
is a singular value decomposition (and an unitary diagonalization) of the scaled kernel matrix $A_W$, $W=diag(w_j)$, where 
$\{w_j\}_{j=1}^N$ are the weigths of a cubature rule $\cubr$.
\end{defin}
\begin{remark}\label{rem:weights}
For what follows it is important to require that $\sum_{j=1}^N w_j = |\Omega|$, which is equivalent to ask that
the cubature rule $\cubr$ is exact at least for the constant functions, $\forall N\in\mathbb{N}$. In the next this 
property will be assumed to hold. 
\end{remark}

\subsection{Preservation of the properties of the continuous basis}
As expected, this basis preserves, in a discretized sense, some interesting properties of the ``natural'' one: 
\begin{prop}\label{pr:prop}
Each weighted svd basis $\U$ has the following properties: 
\begin{enumerate}
 \item $\ u_j(x) = \frac{1}{\sigma_j^2} \sum_{i=1}^N w_i u_j(x_i) \Phi(x,x_i) \approx \frac{1}{\sigma_j^2} \T[u_j](x)$, $\forall\ 1\leqslant j \leqslant N, \ \forall x\in\Omega$
 \item $\U$ is $\Phi$-orthonormal 
 \item $\U$ is $\lw$-orthogonal
 \item $\|u_j \|_{\lw}^2 = \sigma^2_j \quad \forall u_j\in\U$
 \item $\sum_{j=1}^N \sigma_j^2 = \phi(0)\ |\Omega|$
\end{enumerate}
\end{prop}
\begin{pf}
Properties 2 - 3 - 4 can be proved using the gramians as in Prop. \ref{th:gram}:
\begin{eqnarray*}
G_{\U}&=&\Cu^T\cdot A\cdot\Cu = \Cu^T\cdot\Vu =  \Sigma^{-1}  \cdot Q^T \cdot \sqrt{W} \cdot \sqrt{W^{-1}}\cdot Q \cdot \Sigma=I\\
\Gamma_{\U}&=& \Vu^T\cdot W \cdot\Vu = \Sigma \cdot Q^T\cdot \sqrt{W^{-1}} \cdot W\cdot  \sqrt{W^{-1}}\cdot Q \cdot \Sigma = \Sigma^2
\end{eqnarray*}
To prove the first one it suffices to use the definition of $\Cu$ and $\Vu$: indeed from the definition of $\Vu$,
if we denote the $j$-th column of $\Vu$ as ${\Vu}_j$, we get
\begin{eqnarray*}
\Vu &=& \sqrt{W}^{-1} Q \Sigma = \sqrt{W}^{-1} [q_1 \sigma_1,...,q_N\sigma_N ]\\
\Rightarrow&& E_X(u_j) = {\Vu}_j =  \sqrt{W}^{-1} q_j \sigma_j\\
\Rightarrow&& q_j / \sigma_j = 1 / \sigma_j^2 \sqrt{W} E_X(u_j)
\end{eqnarray*}
and the last equality allows to compute each component of $\Cu$ as 
\begin{equation*}
(\Cu)_{i,j}=(\sqrt{W}\cdot Q \cdot \Sigma^{-1})_{i,j} = \sqrt{w_i}\ {\frac{q_j(i)}{\sigma_j}} = {\frac{w_i}{\sigma_j^2}}\ u_j(x_i)  
\end{equation*}
and then the Definition of $\U$ gives 
\begin{eqnarray*}
u_j(x) &=& \sum_{i=1}^N \Phi(x,x_i)\ (\Cu)_{i,j} = \sum_{i=1}^N \Phi(x,x_i)\ {\frac{w_i}{\sigma_j^2}}\ u_j(x_i) =\\
       &=& {\frac{1}{\sigma_j^2}} \sum_{i=1}^N w_i\ \Phi(x,x_i)\ u_j(x_i)
\end{eqnarray*}
where the last term is cleary the approximation of $\T[u_j]$ given by the rule $\cubr$, divided by the corresponding
discrete eigenvalue $\sigma_j^2$.
\\Property 5  is due to a linear-algebra relation: we recall that
\begin{equation*}
\sqrt{W}\cdot A\cdot \sqrt{W} = Q\cdot \Sigma^2 \cdot Q^T 
\end{equation*}
and since the trace of a square matrix is equal to the sum of its eigenvalues, we get
\begin{equation*}
\sum_{j=1}^N \sigma_j^2 = \sum_{j=1}^N w_j\ \Phi(x_j,x_j) = \phi(0)\ \sum_{j=1}^N  w_j  = \phi(0)\ |\Omega|. 
\end{equation*}
This concludes the proof.
\end{pf}

\begin{remark}
In this context, where $\{w_j\}_{j=1}^N$ are cubature weights, the $\lw$-scalar product is a discretization af the 
$L_2(\Omega)$-scalar product. Indeed 
\begin{equation}
(f,g)_{L_2(\Omega)}^2 = \int_{\Omega} f(x) g(x) dx \approx \sum_{j=1}^N w_j f(x_j) g(x_j) =  (f,g)_{\lw}^2 
\end{equation}
and in this sense the property 3 is a discretized version of the corresponding property of the continuous basis.
\\In the same way, property 5 states that
\begin{equation*}
\sum_{j=1}^N \sigma_j^2 = \sum_{j=1}^N w_j\ \Phi(x_j,x_j) = \int_{\Omega} \Phi(x,x)\ dx
\end{equation*}
and exactly the same relation holds for the continuos eigenvalues if $N\rightarrow\infty$, as pointed out in Remark 
\eqref{rem:trace}. 
In this case the integral is exactly approximated by the cubature rule since it is supposed to be exact at least 
for constant functions.  
\end{remark}

\subsection{Completeness}
The basis enjoys another important property: if the data set $X$ if sufficiently big, the basis is complete in the 
native space. Although at the moment it is quite useless in a practical sense, it will be a key property when combined with 
the stability considerations which follows. 
\begin{remark}
A function $f\in\ns$ belongs to $\hs^{\bot}$ if and only if it vanishes on the discrete set $X\subset\Omega$, i.e.
$$\hs^{\bot} = \{ f \in \ns : f(x_i) = 0 \ \forall x_i\in X \}$$
Indeed if $f \in \hs^{\bot}$, i.e. $(f,u_j)_{\Phi}=0\ \forall u_j\in\U$, then using Remark \eqref{rem:intorth}
$$\Int(x) = \sum_{j=1}^N (f,u_j)_{\Phi}\ u_j(x) \equiv 0$$
and $f = \Int = 0$ on $X$.
\\On the other hand, if $f=0$ on $X$, then $\Int = 0$ on $X$, and since $\U$ is a $\Phi$-orthonormal and hence linearly independent set, 
we have $(f,u_j)_{\Phi}=0\ \forall u_j\in\U$.
\end{remark}

\begin{prop}[Completeness]
If the data set $X$ is dense in $\Omega$, the basis $\U$ is a complete basis for $\ns$.
\end{prop}
\begin{pf}
To prove the completeness it suffices to show that if a function $f\in \ns$ is orthogonal to each basis element $u_j\in\U$, 
it is 
the null function. This fact follows immediatly from the previous Proposition, the denseness of $X$ in $\Omega$ and 
the embedding $\ns \hookrightarrow \C(\Omega)$. 
\end{pf}

\section{Interpolation and approximation}
In this section we will describe the interpolation and approximation process based on the weighted svd basis.
This allows to get error bounds and stability estimates that explains the use of such basis.
\subsection{Interpolation}
We recall some basic facts about the interpolation process by kernels. Most of all are the same for each $\Phi$-orthonormal basis,
 as shown in the previous chapter.
\begin{prop}[Interpolant]\label{pr:interpolant}
Let $f\in\ns$ and $X\subset\Omega$, $X=\{x_j\}_{j=1}^N$. Then the interpolant $\Int$ of $f$ on $X$ based on the basis $\U$ 
can be expressed as
\begin{equation}
\Int(x) = \sum_{j=1}^N (f,u_j)_{\Phi} u_j(x)\quad \forall x\in \Omega
\end{equation}
\end{prop}

\begin{prop}[Power function]\label{prop:power}
The Power Function $\Pf$ takes the form
\begin{equation}
\Pf^2= \phi(0) - \sum_{j=1}^N u_j(x)^2  
\end{equation}
\end{prop}
\begin{pf}
The result can be proved using the matrix form of the Power Function given in \eqref{char:pf} for the special case $G_{\U}=I$,
i.e. for an orthonormal basis, but it is more clear to give a direct proof. 
\\By definition, the Power Function is the norm of the error functional
\begin{equation*}
\mathcal{E}_x:\ns \rightarrow \R
\end{equation*}
that maps every function $f\in\ns$ to the interpolation error at a given point $x\in\Omega$, i.e.
\begin{equation*} 
\mathcal{E}_x[f](x) =f(x)-\Int(x) . 
\end{equation*}  
Using the reproducing property of the kernel and the Proposition \eqref{pr:interpolant}, we can rewrite $\mathcal{E}_x$ as
\begin{eqnarray*}
\mathcal{E}_x[f](x) &=& (f,\Phi(\cdot,x))_{\Phi}  - (f, \sum_{j=1}^N u_j(\cdot) u_j(x))_{\Phi} \\
&=&(f,\Phi(\cdot,x)  -  \sum_{j=1}^N u_j(\cdot) u_j(x))_{\Phi}
\end{eqnarray*}
and since the norm of $\mathcal{E}_x$ equals the norm of its Riesz representer, we can conclude that
\begin{eqnarray*}
\Pf^2 &=& \| \Phi(\cdot,x) -  \sum_{j=1}^N u_j(\cdot) u_j(x)\|^2_{\Phi}  \\
&=& \|\Phi(\cdot,x)\|^2_{\Phi} + \sum_{i,j=1}^N u_i(x)u_j(x)( u_i,u_j)_{\Phi} -2 \sum_{j=1}^N u_j(x)( \Phi(\cdot,x),u_j)_{\Phi}\\
&=& \phi(0) - \sum_{j=1}^N u_j(x)^2.  
\end{eqnarray*}
This concludes the proof.
\end{pf}
\begin{remark}
The proof gives also an expansion of the kernel when restricted to act on functions in $\hs$,  
\begin{equation*}
\Phi(x,y)=\sum_{j=1}^N u_j(x) u_j(y) .
\end{equation*}
Indeed $\forall f\in\hs$ we have $f=\sum_{j=1}^N (f,u_j)_{\Phi} u_j(x)$, hence $\forall x\in \Omega$
\begin{equation*}
(f,\Phi(x,\cdot))_{\Phi} = f(x) = \sum_{j=1}^N u_j(x) (f, u_j(\cdot) )_{\Phi} = (f,\sum_{j=1}^N u_j(x) u_j(\cdot) )_{\Phi}   
\end{equation*}
\end{remark}

\subsection{Weighted least-squares approximation}
Another common approach to reconstruct a function from its value at a discrete set $X\subset\Omega$ is to approximate it
in the least-squares sense. 
\\The idea is to use the same sampled data used for the interpolation process, but to relax the interpolation condition
and then to project the function into a subspace of $\hs$, rather than into the full subspace.
\\In this way it is possible to reduce the computational cost of the process, since a smaller basis $\U'\varsubsetneq\U$ is 
involved, and to obtain better results in terms of stability. The hope is to gain these benefits without a serious loss
of convergence speed. 
\\This kind of approximation is also meaningful when the data values are supposed to be affected by noise,
and then an exact recostruction of them makes no sense.
\\In this setting, moreover, the properties of the weighted svd basis provides an additional reason to consider this kind
of approximation: since the eigenvalues of the operator $T_{\Phi}$ decays very rapidly to zero, and the discrete one 
which approximates them are the discrete norms of the basis involved, it is reasonable to consider only the most significant
of them. 
\\\\In fact, we are interested in a weighted least-squares approximation defined as follows:
\begin{defin}\label{def:approx}
Given a function $f\in\ns$, a discrete subset $X\subset\Omega$, a set of cubature weights $\mathcal{W}$ associated with $X$,
a weighted svd basis $\U$ for $\hs$ and a natural number $M\leqslant N = |X|$, the weighted least-squares approximation of order $M$  
of $f$ is the function $\Appr$ that satisfies the codition
\begin{equation*}
\Appr = \argmin_{g\in span\{u_1,\dots,u_M\}} \|f - g\|_{\lw} 
\end{equation*}
\end{defin}
\noindent The reason to use the first basis element is the fact that they are associated to the bigger singular values of the 
scaled kernel matrix $A_W$, and then they allows a more stable and accurate reconstruction of $f$. This relation will be 
stated exactly later when dealing with error bounds.
\\\\In order to compute the approximant, we start to point out a relation between the $\Phi$- and the $\lw$- scalar product, 
that is again a discretized version of a property of the continuous base $\{\varphi_j\}_{j>0}$, stated in Remark 
\eqref{L2Nphi}:    
\begin{prop}
For all $f\in\ns$ and for each $u_j \in \U$, the following relation between the $\Phi$- and the $\lw$-scalar products holds:  
\begin{equation*}
\left(f,u_j\right)_{\Phi} =\frac{1}{\sigma_j^2}\left(f,u_j\right)_{\lw}=\frac{\left(f,u_j\right)_{\lw}}{\left(u_j,u_j\right)_{\lw} } 
\end{equation*}
\end{prop}
\begin{pf}
Using property 1 of Proposition \eqref{pr:prop}, by direct calculations we get the statement:
\begin{eqnarray*}
\left(f,u_j\right)_{\Phi} &=& \left(f, \frac{1}{\sigma_j^2} \sum_{i=1}^N w_i u_j(x_i) \Phi(\cdot,x_i)\right)_{\Phi} =
\frac{1}{\sigma_j^2} \sum_{i=1}^N w_i u_j(x_i) \left( f,\Phi(\cdot,x_i)\right)_{\Phi}\\
&&=\frac{1}{\sigma_j^2} \sum_{i=1}^N w_i u_j(x_i) f(x_i) = \frac{1}{\sigma_j^2\ }\left(f,u_j\right)_{\lw}
\end{eqnarray*}
where $\sigma_j^2 = \left(u_j,u_j\right)_{\lw}$ as shown in the same Proposition.
\end{pf}
\\Now it is possible to express the approximation as a function in $\hs$:
\begin{prop}[Weighted least-squares approximant]
In the notation of Definition \eqref{def:approx}, the approximant is given by
\begin{equation}
\Apprx = \sum_{j=1}^M \frac{(f,u_j)_{\lw}}{\sigma_j^2} u_j(x) = \sum_{j=1}^M (f,u_j)_{\Phi} u_j(x) \quad \forall x\in \Omega
\end{equation}
and then $\Appr$ is nothing else but a truncation of $\Int$.
\end{prop}
\begin{pf}
The second term is simply the orthogonal projection of $f$ into the space generated by $\{u_1,\dots,u_M\}$, since the base
is orthogonal with respect to $(\cdot,\cdot)_{\lw}$ and $\|u_j\|_{\lw}^2=\sigma_j^2$, and it is obviously the element of
 $span\{u_1,\dots,u_M\}$ that minimizes the $\lw$-distance from $f$.
\\The third term is derived as a direct implication of the previous Proposition. 
\end{pf}
\begin{remark}\label{rem:trunc}
We point out that the previous proposition proves that the weighted least-squares approximant $\Appr$ can be obtained 
from the interpolant $\Int$ simply truncating the last $N-M$ coefficients and basis, the ones corresponding to the 
smallest singular values $\sigma_j^2$. This is in opposition to the case of the standard basis of translates, in which 
choosing the bases to neglect corresponds to the choice of a restricted subset $Y\subset X$ where to center the kernel, 
and in general this is a more difficult task. 
\\Moreover, motivated from the results in the following sections, this procedure can be automated when dealing with a 
situation in which very small singular values are expected. In this setting one can leave out the bases which 
correspond to singular values less than a fixed tolerance in order to avoid numerical instability, and then skipping 
automatically from interpolation to least-squares approximation. From a linear algebra point of view, this corresponds
to solve the weighted linear system associated to the interpolation problem using a \textit{total least-squares} method.
\end{remark}
\noindent Using this direct expression for the approximant, it is easy to compute the equivalent of the Power Function for $\Appr$:
\begin{prop}\label{th:lspf}
The norm of the error functional $\mathcal{E}_x^M $ for $\Appr$ takes the form
\begin{equation*}
\left\|\mathcal{E}_x^M\right\|_{\ns^*}^2 =  \phi(0) - \sum_{j=1}^M u_j(x)^2 
\end{equation*}
\end{prop}
\begin{pf}
The proof is the same of the one for the Power Function, where the interpolant is truncated after $M$ terms.
\end{pf}

\subsection{Stability}
In the case of the interpolant, we recall the result of the previous chapter, that holds for each $\Phi$-orthonormal 
basis:
\begin{prop}
Since $\U$ is a $\Phi$-orthonormal basis, it is possible to bound the absolute value of the interpolant as
\begin{eqnarray}
\left| \Intx \right| \leqslant \sqrt{\phi(0)}\ \|f\|_{\Phi} \quad \forall x\in \Omega
\end{eqnarray}
\end{prop}
\noindent The result can be refined for the particular case of an svd-basis:
\begin{prop}\label{pr:stab}
For an svd bases $\U$, the following stability estimates holds:
\begin{equation}
\left| \Intx \right| \leqslant \sqrt{\phi(0)}\ \|f\|_{\Phi}, \quad 
\left| \Apprx \right| \leqslant \sqrt{\phi(0)}\ \|f\|_{\Phi} \quad \forall x\in \Omega
\end{equation}
where in particular
\begin{equation}
\left| \Intx \right| \leqslant \sqrt{\sum_{j=1}^N u_j(x)^2}\ \|f\|_{\Phi}, \quad 
\left| \Apprx \right| \leqslant \sqrt{\sum_{j=1}^M u_j(x)^2}\ \|f\|_{\Phi} \quad \forall x\in \Omega
\end{equation}
\end{prop}
\begin{pf}
It suffices to use the Cauchy-Schwarz inequality for $(\cdot,\cdot)_{\Phi}$ and the $\Phi$-orthonarmality of $\U$:
\begin{eqnarray*}
\left| \Intx \right| &=& \left|\sum_{j=1}^N (f,u_j)_{\Phi} u_j(x)\right| = \left|(f,\sum_{j=1}^N u_j(\cdot)\ u_j(x))_{\Phi}\right|\\
&\leqslant& \left\|\sum_{j=1}^N u_j(\cdot)\ u_j(x)\right\|_{\Phi}\ \| f\|_{\Phi}= \sqrt{\sum_{j=1}^N u_j(x)^2}\ \| f\|_{\Phi}
\end{eqnarray*}
The same for $\Appr$, where the sum stops at $M$. 
\\Now the equality \eqref{prop:power},
\begin{equation*}
\Pf^2=\phi(0) - \sum_{j=1}^N u_j(x)^2  
\end{equation*}
gives exactly the general estimate, since obviously the square of the Power Function is non negative.
\end{pf}

\subsection{Error bounds}
Now we can prove some convergence estimates on the described approximants.
\\It is important to remark that, in the case of the interpolant, there are no difference in the use of a particular
basis, since the spanned subspace $\hs$ in which we project a function $f\in\ns$ clearly does not depend on the
choosen basis. 
\\On the other hand, the fact that we are using this kind of basis allows us to relate the bounds to the continuous
eigenvalues and to their eigenfunctions $\{\varphi_j\}_{j>0}$, which forms a complete basis for $\ns$ and which are related in
a close way to the used kernel $\Phi$. This justifies the choice of sampling the function $f$ on a data set $X$ 
that forms together with a set of weights $\mathcal{W}$ a good cubature rule.
\\Moreover, this remark remains valid in the case of the weighted least-squares approximant, where in addition
the connection between the discrete and the continuos eigenvaues motivates the use of a reduced subspace of $\hs$.
\\\\The first error bound is a simple adaptation of the estimate \eqref{th:bound1} for the case of a 
$\Phi$-orthonormal basis, and then in particular for a weighted svd basis:
\begin{prop}\label{th:bound2}
Let $\Omega\subset\Rn$, let $\Phi\in\mathcal{C}(\Omega\times\Omega)$ be a radial positive definite kernel, and $X\subset\Omega$. 
Then for each $f\in\ns$ and for all $x\in\Omega$, 
\begin{equation}\label{est1}
|f(x) - \Int(x)|^2 \leqslant \left(\phi(0) - \sum_{j=1}^N u_j(x)^2\right)\ \|f\|_{\Phi}^2 
\end{equation}
\end{prop}
\begin{remark}
As shown in Chapter \eqref{ch:general}, in the above estimate the norm of $f$ can be replaced with the norm
of $f-\Int$. In this case, since the basis is $\Phi$-orthonormal, the following equalities holds:
\begin{equation}
\|f -\Int\|_{\Phi}^2 = \|f\|_{\Phi}^2 -\|\Int\|_{\Phi}^2= \|f\|_{\Phi}^2 - \sum_{j=1}^N (f,u_j)_{\Phi}^2 
\end{equation}
The same remains valid also in the next estimates which are derived from this one.
\end{remark}
\noindent\\Using this estimate we can compute a bound on the interpolation and approximation error using the $L_2(\Omega)$-
norm. It is of interest to estimate the reconstruction precision in such norm because it gives not only a punctual, ``worst-case'' bound on the
error, but also a global one.   
\begin{prop}\label{pr:conv1}
Let $\Omega\subset\Rn$ be compact, let $\Phi\in\mathcal{C}(\Omega\times\Omega)$ be a radial positive definite kernel 
and $X\subset\Omega$. Then for each $f\in\ns$ 
\begin{equation*}
\|f - \Int\|_{L_2(\Omega)}^2 \leqslant \left(\left|\Omega\right|\cdot \phi(0)-\sum_{j=1}^N \lambda_j +C_{\Phi,\Omega,X,\mathcal{W}}\cdot\sum_{j=1}^N \|u_j-\varphi_j\|_{L_2(\Omega)}\right) \|f\|_{\Phi}^2 
\end{equation*}
\end{prop}
\begin{pf}
From the embedding Theorem \eqref{th:embed2} we know that $\ns\hookrightarrow L_2(\Omega)$, hence both sides of \eqref{est1}
have finite $L_2(\Omega)$-norm. Thus we can integrate over $\Omega$ the first bound and get
\begin{eqnarray*}
\|f - \Int\|_{L_2(\Omega)}^2 &\leqslant&\mathop{\mathlarger{\mathlarger{\int}}}_{\Omega}\left(\phi(0) - \sum_{j=1}^N u_j(x)^2\right)\ \|f\|_{\Phi}^2 dx\\
&=& \left(|\Omega|\cdot\phi(0)-\sum_{j=1}^N\int_{\Omega}u_j(x)^2 dx\right) \|f\|_{\Phi}^2\\
&=& \left(|\Omega|\cdot\phi(0)-\sum_{j=1}^N\|u_j(x)\|_{L_2(\Omega1)}^2 \right) \|f\|_{\Phi}^2
\end{eqnarray*}
Now we can estimate the $L_2(\Omega)$-norms as follows: using the simple relations
\begin{equation*}
\|\varphi_j\|_{L_2(\Omega)}^2 = \|u_j\|_{L_2(\Omega)}^2 + \|\varphi_j-u_j\|_{L_2(\Omega)}^2 + 2\ (\varphi_j-u_j,u_j)_{L_2(\Omega)}  
\end{equation*}
and
\begin{equation*}
\|u_j\|_{L_2(\Omega)} \leqslant \|\varphi_j-u_j\|_{L_2(\Omega)}+ \|\varphi_j\|_{L_2(\Omega)}
\end{equation*}
we get $\forall j=1,\dots,N$
\begin{eqnarray*}
-\|u_j\|_{L_2(\Omega)}^2&=&-\|\varphi_j\|_{L_2(\Omega)}^2+\|\varphi_j-u_j\|_{L_2(\Omega)}^2 +2\ (\varphi_j-u_j,u_j)_{L_2(\Omega)}\\
&\leqslant&-\|\varphi_j\|_{L_2(\Omega)}^2+\|\varphi_j-u_j\|_{L_2(\Omega)}^2 +2\ \|\varphi_j-u_j\|_{L_2(\Omega)}\ \|u_j\|_{L_2(\Omega)}\\
&\leqslant&-\|\varphi_j\|_{L_2(\Omega)}^2+\|\varphi_j-u_j\|_{L_2(\Omega)}\left(3\ \|\varphi_j-u_j\|_{L_2(\Omega)}+\|\varphi_j\|_{L_2(\Omega)}\right)
\end{eqnarray*}
and from \eqref{theOperator} we know that $\|\varphi\|_{L_2(\Omega)}^2=\lambda_j$. To conclude it suffices to bound in an obvious way the right-hand 
side term as 
\begin{equation*}
3\ \|\varphi_j-u_j\|_{L_2(\Omega)}+\|\varphi_j\|_{L_2(\Omega)}\leqslant 3\left(\max_{j=1,\dots,N} \|\varphi_j-u_j\|_{L_2(\Omega)}\right) + \sqrt{\lambda_1}=:C_{\Phi,\Omega,X,\mathcal{W}}
\end{equation*}
since the eigenvalues $\{\lambda_j\}_{j>0}$ are non increasing.
\end{pf}
\newline\noindent The same estimate remains valid for the approximant $\Appr$ if $N$ is replaced by $M$, as a consequence of the 
Proposition \eqref{th:lspf}:
\begin{prop}\label{pr:errorls}
Let $\Omega\subset\Rn$ be compact, let $\Phi\in\mathcal{C}(\Omega\times\Omega)$ be a radial positive definite kernel, 
$X\subset\Omega$, $|X|=N\in\mathbb{N}$ and $M\leqslant N$. Then for each $f\in\ns$ 
\begin{equation*}
\|f - \Appr\|_{L_2(\Omega)}^2 \leqslant \left(\left|\Omega\right|\cdot \phi(0)-\sum_{j=1}^M \lambda_j +C_{\Phi,\Omega,X,\mathcal{W}}\cdot\sum_{j=1}^M \|u_j-\varphi_j\|_{L_2(\Omega)}\right) \|f\|_{\Phi}^2 
\end{equation*}
\end{prop}
\noindent We point out that these estimates involve two terms: the first one is 
\begin{equation*}
|\Omega|\cdot \phi(0)-\sum_{j=1}^N \lambda_j
\end{equation*}
and it is related only on the kernel, the domain and the dimension $N\in\mathbb{N}$ of the approximation subspace $\hs$.
From the Remark \eqref{rem:trace} on the Theorem \eqref{th:mercer} we know that for $N\rightarrow\infty$ the above term
vanishes, and moreover the eigenvalues are positive and orderer in a decreasing way. Hence this term measures how the 
truncated series approximates the full one, or in other words how the degenerate kernel
\begin{equation*}
\sum_{j=1}^N \lambda_j \varphi_j(x) \varphi_j(y)  ,\ x,y\in\Omega
\end{equation*}
approximates the original kernel $\Phi(x,y)$.
\\The second term is
\begin{equation*}
C_{\Phi,\Omega,X,\mathcal{W}}\cdot\sum_{j=1}^N \|u_j-\varphi_j\|_{L_2(\Omega)}
\end{equation*}
and it depends also on the cubature rule $\cubr$. It measures the convergence rate of the Nyst\"om method based on the
rule, and gives informations on how well the discrete basis $\U$ approximates the continuous one.
\begin{remark}
It can be useful to refer to the cubature error in terms of the $\|\cdot\|_{\infty}$-norm. It can be done in an obvious way
since $\Omega$ is compact and then
\begin{equation*}
\sum_{j=1}^N \|u_j-\varphi_j\|_{L_2(\Omega)} \leqslant \sqrt{|\Omega|}\ \sum_{j=1}^N \|u_j-\varphi_j \|_{\infty} 
\end{equation*}
\end{remark}

\noindent For the weighted least-squares approximant, it make sense to consider also another type of error measurement. Indeed in this 
case the data-sites set $X\subset\Omega$ is not used to interpolate the function $f$, but works as a sample set, so the
pointwise distance between $f$ and $\Appr$ on $X$ is not zero. We can bound this quantity as shown in the next Proposition:
\begin{prop}
Let $\Omega\subset\Rn$, let $\Phi\in\mathcal{C}(\Omega\times\Omega)$ be a radial positive definite kernel, 
$X\subset\Omega$, $|X|=N$, and $M<N$. Then for each  $f\in\ns$
\begin{equation}
\|f - \Appr\|_{\lw} \leqslant \left(\sum_{j=M+1}^N \sigma_j^2 \right)^{\frac{1}{2}} \|f \|_{\Phi}
\end{equation}
\end{prop}
\begin{pf}
We start again from the bound in Proposition \eqref{th:bound2}: in this case the finiteness of the $\lw$-norm of both sides is obvious,
since every function involved in the estimate is continuous on $\Omega$.
Acting as in the previous Theorem we can evaluate both sides in $X$ and sum the obtained values, weighted with the 
weights $\mathcal{W}$.
We get
\begin{eqnarray*}
\|f - \Appr\|_{\lw}^2 &\leqslant& \ \|f\|_{\Phi}^2  \sum_{i=1}^N w_i\ \left(\phi(0) - \sum_{j=1}^M u_j(x_i)^2\right)\\
&=& \|f\|_{\Phi}^2 \left( \phi(0) \sum_{i=1}^N w_i  - \sum_{j=1}^M  \sum_{i=1}^Nw_i\ u_j(x_i)^2\right)\\
&=& \|f\|_{\Phi}^2 \left( \phi(0) \sum_{i=1}^N w_i  - \sum_{j=1}^M  \| u_j\|_{\lw}^2\right)\\
&=& \|f\|_{\Phi}^2 \left( \phi(0) |\Omega|  - \sum_{j=1}^M   \sigma_j^2\right)
\end{eqnarray*}
Where we used property 4 of Proposition \eqref{pr:prop} to compute the $\lw$-norm of the basis functions and Remark 
\eqref{rem:weights} to compute the sum of the weights.
To get the desired bound we recall that
\begin{equation*}
\sum_{j=1}^N \sigma_j^2 = \phi(0)\ |\Omega| 
\end{equation*}
again as stated in Proposition \eqref{pr:prop}.
\end{pf}
\begin{remark}
 The last result can be interpreted also in another way. In fact, it gives a bound on how the weighted least-squares approximant and
the interpolant differs on the data-sites set $X$. Indeed, since $f(x_i) = \Int(x_i)$ $\forall x_i\in X$, we get 
$\forall f \in \ns$
\begin{equation}
\|\Int - \Appr\|_{\lw} \leqslant \left(\sum_{j=M+1}^N \sigma_j^2 \right)^{\frac{1}{2}} \|f \|_{\Phi}
\end{equation}
Clearly this estimate doesn't give informations on the distance between the two approximants on the set 
$\Omega\setminus X$. However this quantity can be computed as in the above estimates:
\begin{eqnarray*}
|\Intx-\Apprx|&=&\left|\sum_{j=1}^N (f,u_j)_{\Phi} u_j(x)- \sum_{j=1}^M (f,u_j)_{\Phi} u_j(x)\right|\\
&=& \left|\sum_{j=M+1}^N (f,u_j)_{\Phi} u_j(x)\right|\\
&\leqslant& \|f\|_{\Phi}\left( \sum_{j=M+1}^N u_j(x)^2\right)
\end{eqnarray*}
and the partial sum can be bounded as in the previous Theorems.
\end{remark}
\begin{remark}
The trade-off principle between stability and convergence explained in Chapter \eqref{ch:general} can be viewed in this context 
as follows: we have
\begin{eqnarray*}
|\Intx|^2&\leqslant&\left(\sum_{j=1}^N u_j(x)^2\right)\ \|f\|_{\phi}^2\\
|\Intx-f(x)|^2&\leqslant&\left(\phi(0)-\sum_{j=1}^N u_j(x)^2\right)\ \|f\|_{\phi}^2
\end{eqnarray*}
and the same for $\Appr$ if $N$ is replaced by $M<N$.
\\Hence for convergent approximant, namely for approximant for which the power function converges to zero, we have necessarily
\begin{equation*}
\sum_{j=1}^N u_j(x)^2 \rightarrow \phi(0)
\end{equation*}
that is, the constants in the stability bounds in Proposition \eqref{pr:stab} are maximized. 
\end{remark}

\chapter{Numerical experiments}\label{ch:numerics}
In this chapter we will present some numerical experiments which show the actual behaviour of our basis.
\\We will point out different features that can be relevant in the choice of a method, 
and in particular the ones concerning stability and convergence speed. 
\\The following tests take place in the setting described in Section \eqref{num:general}.
\\\\The code is mainly written in \Mat, using in some parts the software present in the book \cite{Fass}. The most 
performance-critical parts are written in \cpp, using the \Mat \textit{MEX} interface \cite{Mex} and the 
linear-algebra library \textit{Eigen} \cite{Eigen}. 

\section{General setting}\label{num:general}
The approximant strictly depends on the set $\Omega\subset\Rn$, on the kernel 
$\Phi\in\mathcal{C}(\Omega\times\Omega)$, on the data-sites set $X\subset\Omega$ and on the function $f:\Omega\rightarrow \R$
that we try to reconstruct. In this section we describe the general choices made on this elements for our experiments.  
 
\subsection{Approximation domain}\label{sect:domains}
\begin{figure}
\center
\includegraphics[width=6cm]{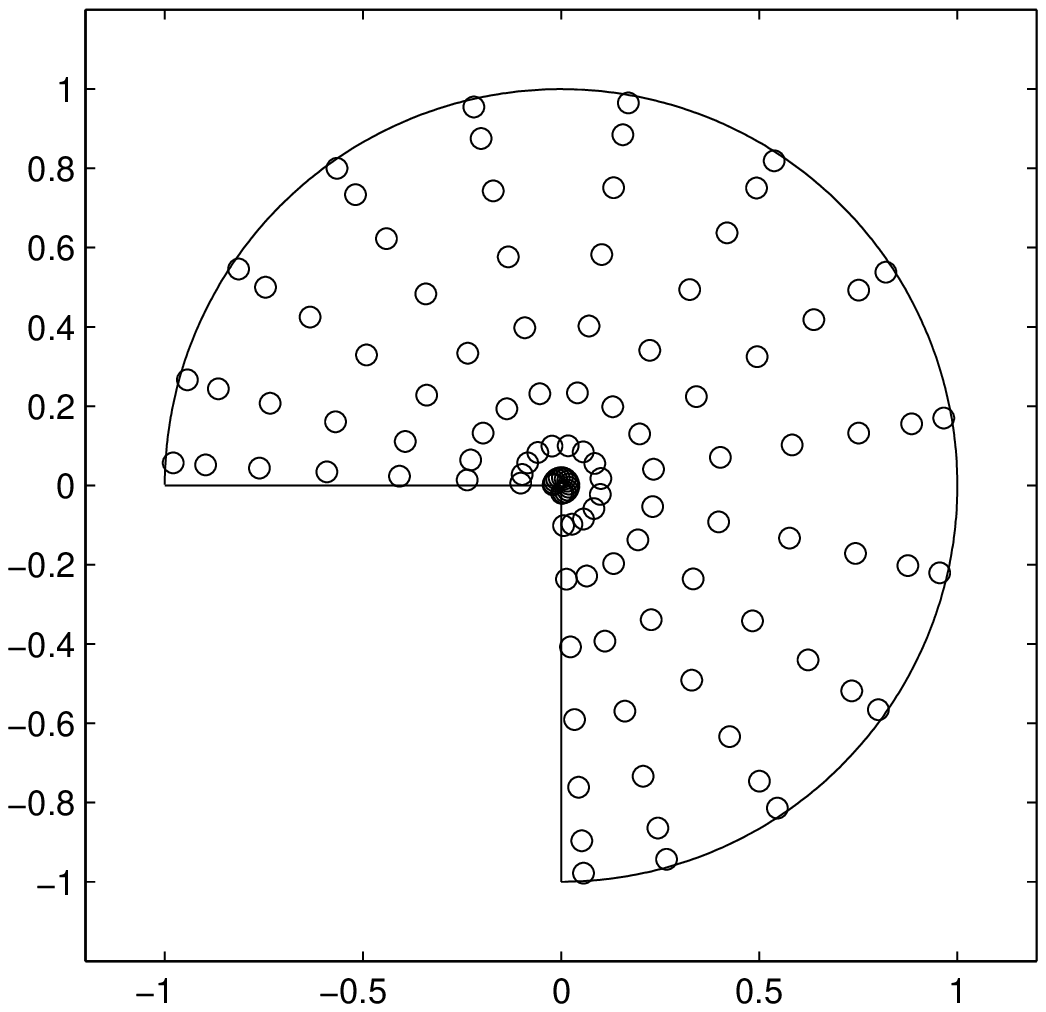}
\includegraphics[width=6cm]{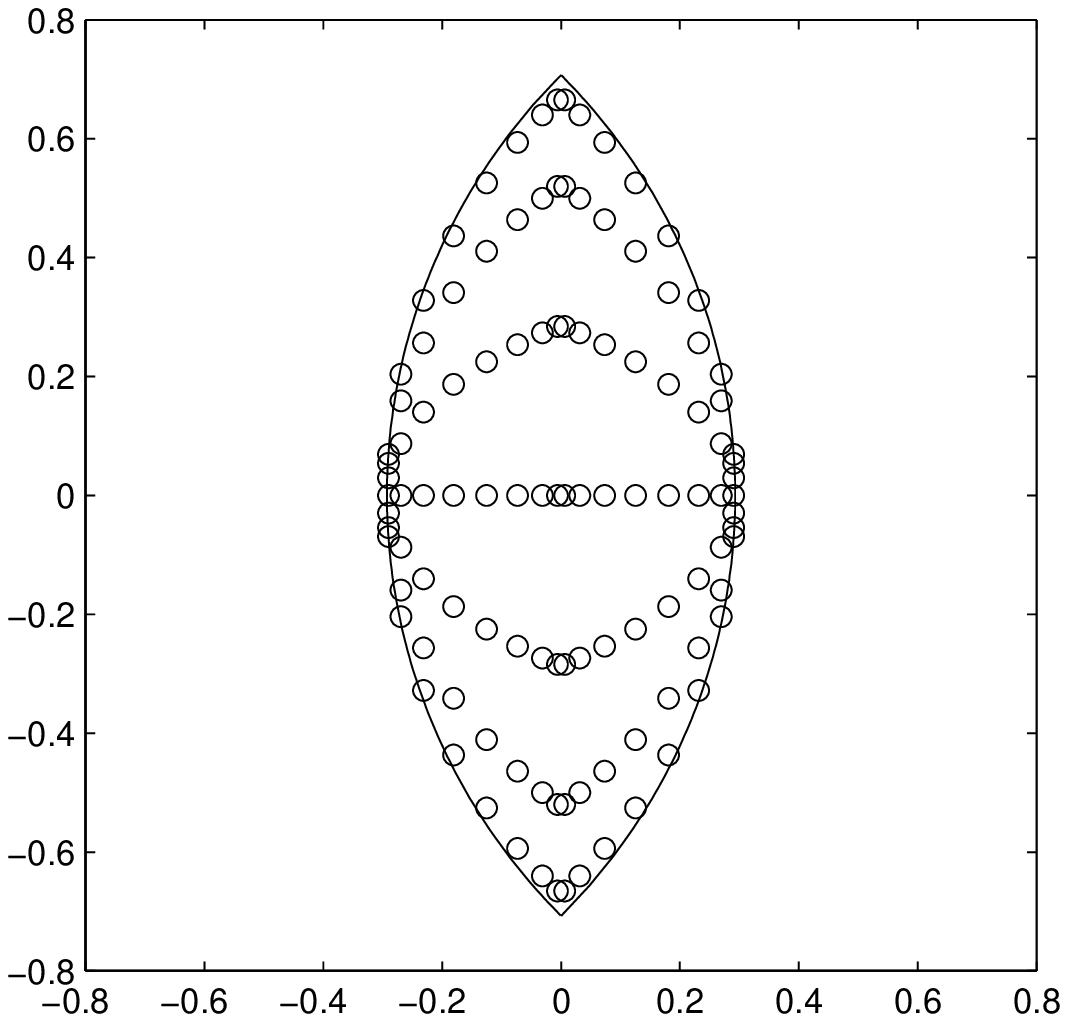}
\caption{The cutted disk and the lens as defined in the Section \eqref{sect:domains}, with an example of trigonometric gaussian
points, as defined in Section \eqref{sect:dataS}}
\label{fig:cuttedandlens}
\end{figure}

The sets $\Omega\subset\R^2$ used are the following:
\begin{itemize}
 \item The square $\Omega_1=[0,1]\times[0,1]$ 
 \item The disk $\Omega_2$ with center $C=\left(\frac{1}{2},\frac{1}{2}\right)$ and radius $R=\frac{1}{2}$ 
 \item The \textit{cutted disk} $\Omega_3$, i.e. the unit disk centered in zero with the third quadrant cutted away.
 It is the domain depicted on the left in Figure \ref{fig:cuttedandlens}.
 \item The \textit{lens} $\Omega_4$ defined as the intersection of two disk with centers  $C=\left(-\frac{\sqrt{2}}{2},0\right)$ and
$c=\left(\frac{\sqrt{2}}{2},0\right)$ and radii $R=r=1$. It is the domain depicted on the right in Figure \ref{fig:cuttedandlens}. 
\end{itemize}

\subsection{Data sites}\label{sect:dataS}
For the standard basis of translates we use two distributions of points: equally spaced points and Halton points. The 
Halton points are used since they are well distributed without being equally spaced.  
\\For our basis we need a data-sites set $X\subset\Omega$ such that $\cubr$ is a cubature rule. We use product Gauss-Legendre 
points for the square $\Omega_1$, and \textit{trigonometric gaussian} points on the other sets in $\R^2$. 
\\The latter set of cubature points was recently presented in the papers \cite{ViDF,ViBo}, and can be obtained for a wide
class of domains, much more general than the one used here. We use them because they provide a high-accuracy cubature 
rule while being sufficiently uniform in $\Omega$. {\Mat} functions to compute this points can be found in the site 
\cite{Vianello}.    
\subsection{Test functions}
The functions we try to reconstruct are typical test functions in the context of approximation theory:
\begin{itemize}
\item the bivariate Franke function $f_F$
\item an oscillatory function $f_o(x,y)=\cos(20 (x+ y))$
\item a function with a derivative discontinuity at $x=y$, $f_s(x,y) = e^{\|x-y\|_2}-1 $
\item a function $f_{\mathcal{N}}$ belonging to the native space of the gaussian kernel, obtained as a linear combination of the kernel 
centered on some points in $\Omega$, for a fixed shape parameter. We use them to test the behaviour of the 
approximant for functions in the native space. 
\end{itemize}

\subsection{Kernels}
We use three different kernels among the ones described in Table \eqref{tb:kernels}.
\\The choice is motivated from the different behaviour of the eigenvalues $\{\lambda_j\}_{j>0}$ of the integral operator $T_{\Phi}$ associated with these basis functions. Indeed, 
altough we know from Theorem \eqref{th:mercer} that the continuous eigenvalues accumulate to zero, the speed in which they 
decay is clearly not the same for different kernels.  
\\The basis functions are the gaussian (fast decay to zero), the IMQ (slower decay) and the 3MAT (slow decay). Nevertheless,
we point out that also the choice of the shape parameter $\varepsilon$ strongly influences this speed.
\\We can expect that this difference reflects on the approximation and on its stability.   
\begin{figure}
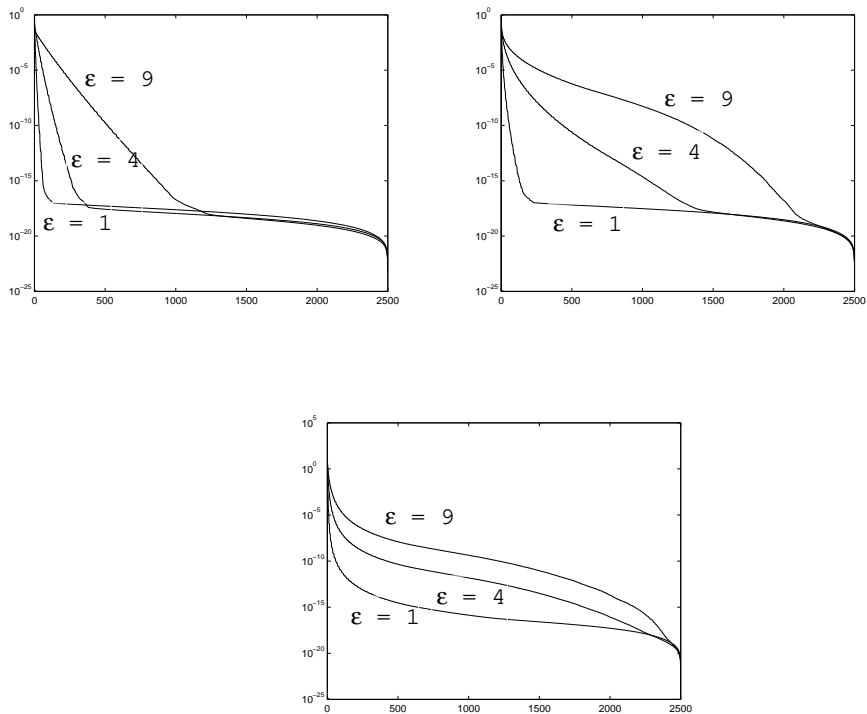

\includegraphics[width=6cm]{estimatedEig.disk.gauss.1.4.9}
\includegraphics[width=6cm]{estimatedEig.disk.imq.1.4.9}\\
\center
\includegraphics[width=6cm]{estimatedEig.disk.mat3.1.4.9}
\caption{Eigenvalues for the gaussian (top left), the IMQ (top right) and the 3MAT (bottom) kernels for shape parameter $\varepsilon = 1, 4, 9$, 
computed on the disk $\Omega_2$ using the Nystr\"om method based on $50^2$ trigonometric Gauss-Legendre cubature points 
and weights.  }
\label{fig:eigen}
\end{figure}

\section{Comparison between the interpolant and the weighted least-squares approximant}
\begin{figure}
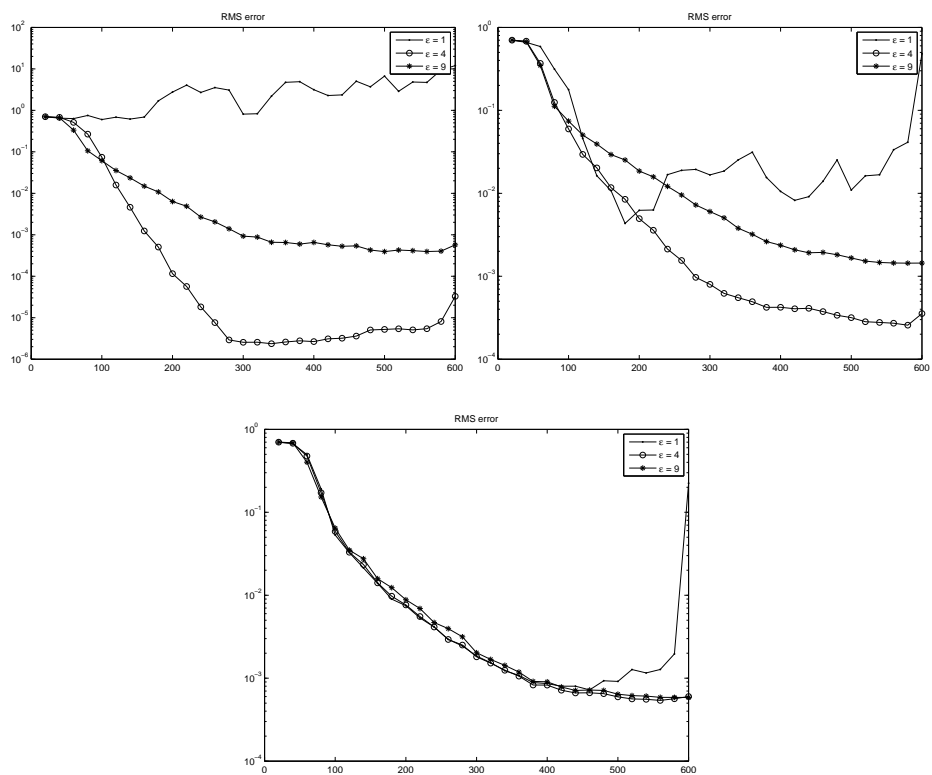

\center
\includegraphics[width=6cm]{IntVsApp.disk.gauss.1.4.9..nerr}
\includegraphics[width=6cm]{IntVsApp.disk.imq.1.4.9..nerr}\\
\center
\includegraphics[width=6cm]{IntVsApp.disk.mat3.1.4.9..nerr}
\caption{RMS error for the reconstruction of $f_o$ on $\Omega_2$ using $\Appr$ for different Ms and different shape parameters, using 
the gaussian kernel (top left), the IMQ (top right) and the 3MAT kernel (bottom).}
\label{fig:IntVsAppr}
\end{figure}

As pointed out in the Remark \eqref{rem:trunc}, we know that we can compute the weighted least-squares approximant as a 
truncation of the interpolant.
\\This reduction increases the error residual as shown in Proposition \eqref{pr:errorls}, but in the cases in which the 
smaller eigenvalues are under a certain tolerance, we can expect that a truncation does not affect too much the 
approximation capability. Furthermore, altough Proposition \eqref{pr:stab} proves the stability of our basis, in some 
limit situations we can expect that the influence of the smallest eigenvalues produces numerical instability that cannot
 be completely controlled. 
\\In the next example we compare the approximation error produced using the full interpolant and some reduced weighted least-squares
approximant. We reconstruct the oscillatory function $f_o$ on the disk $\Omega_2$ using the three kernels with 
$\varepsilon=1, 4, 9$, starting from $600$ trigonometric gaussian centers, and then truncating the basis for $M \in \{0,20,\dots,600\}$.    
\\To measure the accuracy of the reproduction obtained with this process, we compute the root-mean-square errors (RMSE)
on a uniform grid. Figure \eqref{fig:IntVsAppr} shows the results obtained.
\\\\The results reflect the consideration on the eigenvalues made in the last section. Indeed, we can see that for the 
3MAT kernel the interpolant remains stable for each $\varepsilon$, the last iterations for $\varepsilon = 1$ apart,
 the IMQ becomes instable for $\varepsilon = 1$, while the gaussian presents some instability also for $\varepsilon = 4$.
In the instable cases, there is a clear gain using a truncated version of the interpolant, that is a least-squares approximant $\Appr$ for some $M$. Table
\eqref{tb:M} shows the index $M$ such that $\Appr$ provides the best approximation of $f_o$.
\begin{table}
\centering
\begin{tabular}{||l|c|c|c||}
\hline 
&$\varepsilon=1$&$\varepsilon=4$&$\varepsilon=9$\\
\hline
Gaussian&100&340&500\\
IMQ&180&580&580\\
3MAT&460&560&580\\
\hline
\end{tabular} 
\label{tb:M}
\caption{Optimal $Ms$ for different kernels and shape parameter, i.e. indexes such that the weighted least-squares 
approximant $\Appr$ provides the best approximation of the test function $f_o$ on $\Omega_2$.}
\end{table}
\\\\A special situation occurs for the gaussian with $\varepsilon = 1$, where the reducing process does not suffices
to avoid instability, as expected from the distribution of the eigenvalues shown in top-left Figure \eqref{fig:eigen}:
for this parameter the eigenvalues are almost all under the machine precision. Moreover, for $\varepsilon=1$ the gaussian
becomes too flat, and then there is no hope to reconstruct an oscillatory function.

\section{Comparison with the standard basis}
In the following tests we compare the approximations obtained from the standard basis of translates and from our basis. 
\\Since acoording to the Section \eqref{sec:stdconvstab} we know that the stability of the standard basis is strictly related to the shape 
parameter $\varepsilon$, at first we compare the two methods for different fixed values of $\varepsilon$ and different 
kernels, considering situations in which the standard interpolant becomes seriously instable as well as more stable 
cases. Then we repeat the same tests for an optimized shape parameter $\varepsilon^*$. 
\subsection{Fixed shape parameter}
\begin{figure}
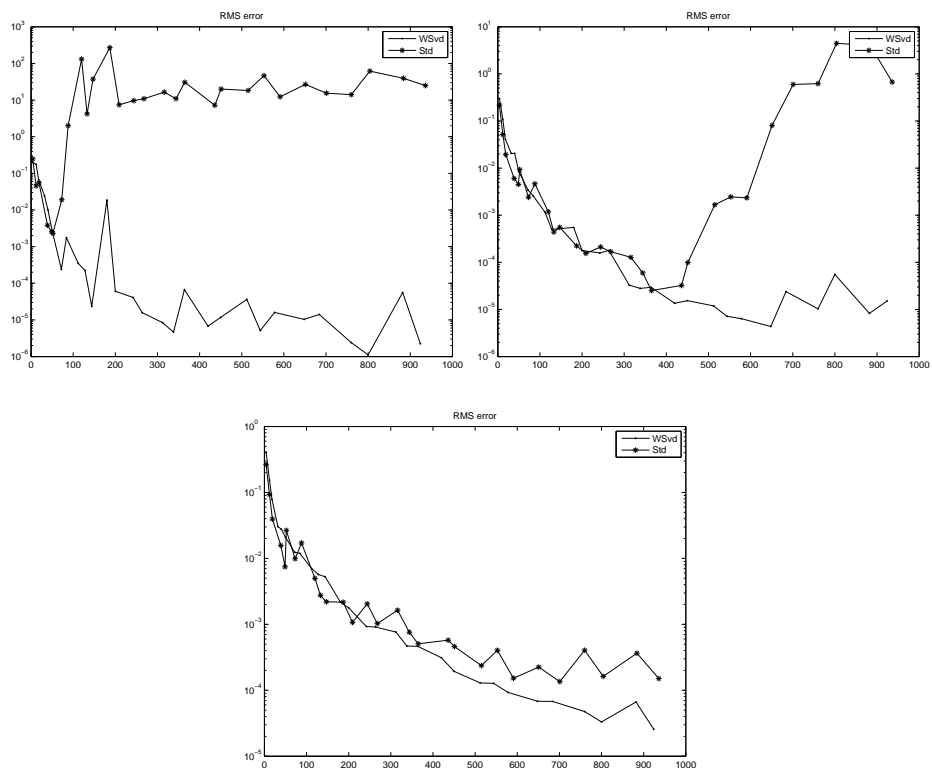

\center
\includegraphics[width=6cm]{WSvdVsStd.lens.u.imq.1.0.franke.2.30.nerr}
\includegraphics[width=6cm]{WSvdVsStd.lens.u.imq.4.0.franke.2.30.nerr}\\
\center
\includegraphics[width=6cm]{WSvdVsStd.lens.u.imq.9.0.franke.2.30.nerr}
\caption{RMS error for the reconstruction of $f_F$ on the lens $\Omega_4$ using the IMQ kernel with the standard basis and our basis.
Test for different shape parameter: $\varepsilon=1$ (top left), $\varepsilon=4$ (top right) and $\varepsilon=9$ (bottom).}
\label{fig:wSvdVsStd}
\end{figure}

\begin{figure}
\center
\includegraphics[width=6cm]{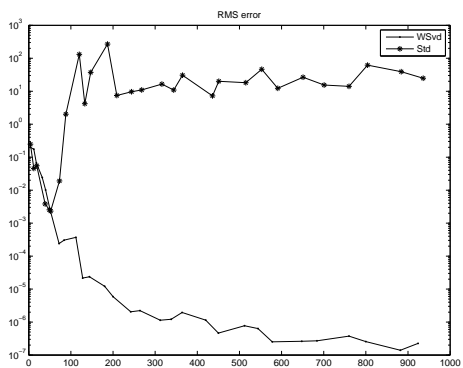}
\caption{RMS errors for the reconstruction of $f_F$ on the lens $\Omega_4$ using the IMQ kernel and $\varepsilon=1$, using 
the interpolant based on the standard basis and the weighted least-squares approximant $\Appr$ with $M$ such that 
$\sigma_M<10^{-17}$}
\label{fig:wSvdVsStdtrunc}
\end{figure}

In this example we try to reconstruct the Franke function $f_F$ on the lens $\Omega_4$ using the IMQ kernel.
\\The test compares the results obtained from the interpolant based on the standard basis centered on an uniform grid and
the one based on our basis, centered on a trigonometric gauss set. The reconstruction is repeated for $\varepsilon=1, 4, 9$
and for data sites sets $X_N\subset \Omega_4$, with $N=|X_N|<1000$. 
\\The RMS errors are reported in Figure \eqref{fig:wSvdVsStd}.
\\We can see that in the stable case, namely for $\varepsilon=9$, there is only a small difference between the two basis, 
altough for $N>500$ the standard interpolant does not gain accuracy. 
\\For $\varepsilon = 1, 4$, altough for small data sets $X_N$ the two basis does not behave much different, 
 when $N$ becomes too big the standard basis becomes instable and a growing of the data-sites set does not lead
to a more accurate reconstruction. On the other hand, the interpolant based on our basis presents a convergent behaviour 
for each shape parameter, even if it is also clearly influenced in the rate of convergence. 
\\This feature can be useful since, at least in the considered cases, there is no need to choose a particular 
$\varepsilon$ to guarantee convergence, even if slow.
\\\\Furthermore, when a small shape parameter influences too much the stability of the interpolant, we can use the reduced
weighted last-squares approximant $\Appr$, as discussed in the previous section. The approximation process for 
$\varepsilon=1$ is repeated using $\Appr$ instead of $\Int$, with $M$ such that $\sigma_M<10^{-17}$. The result is shown
in Figure \eqref{fig:wSvdVsStdtrunc}. The aproximant is clearly more stable, while the convergence rate is not reduced.
\subsection{Optimized shape parameter}\label{sec:optim}
\begin{figure}
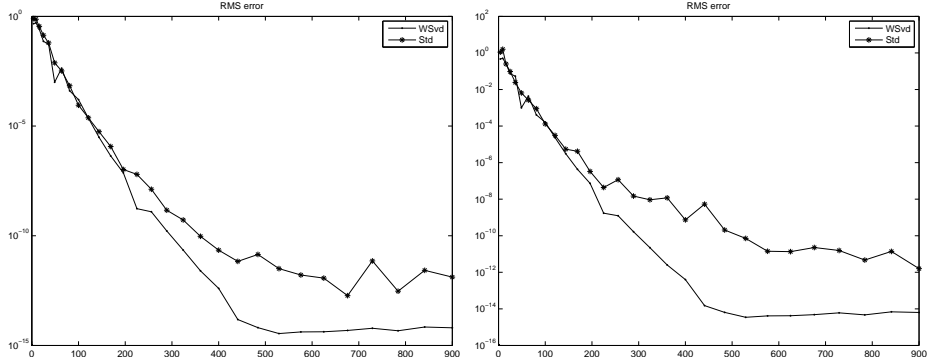

\center
\includegraphics[width=6cm]{WSvdVsStdOpt.square.u.gauss.4.1e-17.exp3.2.30.u.nerr}
\includegraphics[width=6cm]{WSvdVsStdOpt.square.u.gauss.4.1e-17.exp3.2.30.h.nerr}\\
\caption{RMS error for the reconstruction of $f_{\mathcal{N}}$ on the square $\Omega_1$ using the gaussian kernel with the 
standard basis and an optimized shape parameter $\varepsilon*$ and our basis with $\varepsilon=4$. The standard interpolant
is computed using equally spaced points (on the left) and Halton's points (on the right). Our basis is truncated at $M$ such
that $\sigma_M<10^{-17}$}
\label{fig:wSvdVsStdnative}
\end{figure}

A possible solution for the instability of the standar basis is to optimize the shape parameter. In practice, for 
a fixed point distribution, a fixed kernel and a fixed test function, one tries to find the parameter that minimizes the residual error.
\\In the following examples we realize this optimization using the so-called \textit{leave one out} method. The idea
is to compute the interpolant $\Int$ on the full set $X\subset\Omega$ and the $N$ interpolants $P[f]_i$ on the reduced sets 
$X_i=X\setminus\{x_i\}\ \forall\ i\in\{1,\dots,N\}$, for different shape parameter $\varepsilon\in E$, $E\subset\R$, and then to choose the 
optimal $\varepsilon^*$ defined as
\begin{equation*}
\varepsilon^*=\argmin_{\varepsilon\in E} \max_{i=1,\dots,N} \left|\Int(x_i)-P[f]_i(x_i)\right|
\end{equation*}
\\We remark that this optimization is quite expensive in terms of computational time, and cannot be 
performed once for all, but has to be repeated if the data-sites set increases. Moreover, there are cases in which a 
particular choice of $\varepsilon$ is motivated by theoretical reason.
\\To examine this situation, we use as a test function an element of the native space of the gaussian  
$\Phi_4 (x,y) := \exp(-4^2 \|x-y\|_2^2) $ on the square $\Omega_1$, i.e. the function 
\begin{equation*}
f_{\mathcal{N}}(x) = -2 \Phi_4(x,(0.5,0.5))+\Phi_4(x,(0,0))+3 \Phi_4(x,(0.7,0.7))\ \forall x\in[0,1]^2 
\end{equation*}
The RMS errors are plotted in Figure \eqref{fig:wSvdVsStdnative}, using uniform points (on the left) and Halton points
 (on the right) as centers of the standard basis. It is clear that a good choice of the shape parameter reduces 
the instability of the standard interpolant, altough it does not suffices to avoid it completely. On the other hand, 
the stability of our basis, together with the truncation at $M$ such that $\sigma_M<10^{-17}$, allows to use the ``right'' shape parameter for each number of centers, and this leads to an 
approximant that converges to the sampled function with a tolerance near to the machine precision. 
Table \eqref{tb:wSvdVsStdnative} shows the RMS errors for different numbers of data sites, together with the optimal 
parameter $\varepsilon^*$ selected by the leave-one-out optimization.
\begin{table}
\centering
\begin{tabular}{||c|c|c|c|c|c||}
\hline
N&196&324&529&729&900\\
\hline 
Std - e&$1.05\cdot 10^{-7}$&$5.23\cdot10^{-10}$&$3.17\cdot10^{-12}$&$7.15\cdot10^{-12}$&$1.30\cdot10^{-12}$\\
$\varepsilon^*$&$3.75$&$3.81$&$3.84$&$3.87$&$3.98$\\
\hline
Std - H&$3.30\cdot10^{-7}$&$9.31\cdot10^{-9}$&$7.12\cdot10^{-11}$&$1.56\cdot10^{-11}$&$1.59\cdot10^{-12}$\\
$\varepsilon^*$&$3.75$&$3.84$&$3.89$&$3.92$&$3.95$\\
\hline
W-Svd&$7.37\cdot10^{-8}$&$2.23\cdot10^{-11}$&$3.48\cdot10^{-15}$&$6.08\cdot10^{-15}$&$6.37\cdot10^{-15}$\\
\hline
\end{tabular} 
\label{tb:wSvdVsStdnative}
\caption{RMS errors for the approximation described in Section \eqref{sec:optim} obtained using our basis (W-Svd), 
the standard basis centered on equally spaced points (Std - e) and on the Halton points (Std - H), together with 
the optimal shape parameter used for the standard basis. Values for $N=|X|$ as in the first row.}
\end{table}

\section{Comparison with the Newton basis}\label{sec:newton}
\begin{figure}
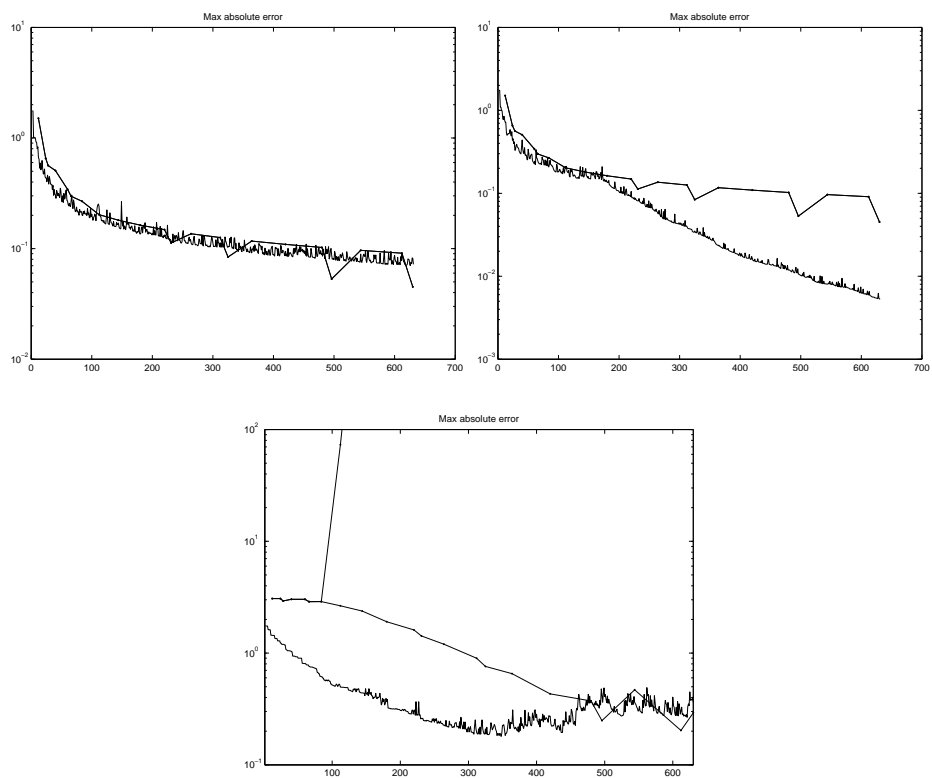

\center
\includegraphics[width=6cm]{WSvdVsNewton.shab.u.wen21.0.5.0.shab.3.25.h.0.01.merr}
\includegraphics[width=6cm]{WSvdVsNewton.shab.u.wen21.0.5.0.shab.3.25.h.0.05.merr}\\
\center
\includegraphics[width=6cm]{WSvdVsNewton.shab.gauss.9.1e-10.shab.3.25.merr}
\caption{Absolute error for the reconstruction of $f_{s}$ on the domain $\Omega_3$ using the Newton basis and our basis 
(dotted lines).
The kernel used are the Wendland's W21 with $\varepsilon = 2$ (first row), using different starting grid for the Newton 
basis, and the gaussian with $\varepsilon = 4$ in the second row (the error for the full interpolant is cutted from the 
plot after it becomes bigger than $10^2$).}
\label{fig:New}
\end{figure}

In the paper \cite{SchPa} the general change of basis described in the Chapter \eqref{ch:change} was the starting point
to create a \textit{Newton basis} for the native space $\ns$. Among other properties partially enjoied also by our 
basis, the Newton basis can be computed recursively, i.e. if we add to $X$ a further sample point $x_{N+1}\notin X$ it suffices
to compute the basis element corresponding to the new point. Moreover, an adaptive point-selection algorithm is provided 
to choose the point to add using also the information given by the sampled values of the test function.
\\{\Mat} programs for the adaptive calculation of the interpolant based on this basis can be found in \cite{SchSite}.
\\\\To compare the two bases we try to reconstruct the function $f_s$ on the cutted disk $\Omega_3$ using the gaussian kernel
and a Wendland kernel W21. We choose to use the latter since it is already present in the mentioned {\Mat} functions and 
because it provides interesting results which are explained in the next.
\\In both cases the adaptive algorithm is able to detect the derivative discontinuity at $x=y$, and to concentrate the data sites near this area.
\\\\The test for the first kernel with $\varepsilon = 4$ shows a much better behaviour of the Newton basis. Indeed, a small 
number of points suffices to introduce high instability in our basis. This can be a consequence of the distribution of 
the data-sites depicted in Figure \eqref{fig:cuttedandlens}, which are too concentrated in zero to produce a well-conditioned kernel matrix.
To avoid this, we repeat the experiment using a weighted least-squares approximant $\Appr$ with the $M$ that provides the
best approximation, namely $M$ such that $\sigma_M<10^{-10}$, and in this case the residual decreases as $N$ grows. The 
maximal absolute erros for $N\leqslant 625$ are depicted in the bottom of Figure \eqref{fig:New}.  
\\\\The second kernel is used with $\varepsilon = 2$, and in this case our basis and the Newton basis provides an almost
equivalent decrease of the maximal absolute error in the range under consideration, as reported in the top left of 
Figure \eqref{fig:New}. 
\\Here the Newton basis behaves in a quite unexpected way, since if the uniform grid from which the data-sites are selected by the algorithm is reduced, namely
if the one dimensional grid size changes from $\Delta_x=0.01$ to $\Delta_x=0.05$, the accuracy of the interpolant becomes 
strictly better, as depicted in the top right part of Figure \eqref{fig:New}.

\chapter{Conclusions and further work} 
This thesis presents a way to construct a new kind of stable basis for the RBF approximation.
\\The basis is derived from the procedure described in Chapter \eqref{ch:change}, 
and is build to have useful properties related to its $\Phi$-orthonormality.
\\Furthermore, the particular approach used here connects this discrete basis with the ``natural'' one described in
the Theorem \eqref{th:mercer}, and this allows to relate some functional property of the kernel to the approximant itself.
\\\\In this setting, a more deep study could lead to a stronger use of the information provided by the kernel and the domain.
In particular  the convergence estimate of Proposition \eqref{pr:conv1} can be refined considering the rate of convergence
 to zero of the eigenvalues of the operator $T_{\Phi}$ and the property and the convergence rate of the Nystr\"om method
based on the setting of the problem, namely the choosen cubature rule, the kernel $\Phi$, the shape parameter $\varepsilon$
and the set $\Omega$.    
\\\\As regards stability, the experiments presented in Chapter \eqref{ch:numerics} confirms the good behaviour expected from
Proposition \eqref{pr:stab}. In particular our basis allows to treat approximants based on a relatively big number of 
points also for not optimized shape parameters and on quite general sets. This feature can be enforced thanks to the 
possibility to compute a weighted least-squares approximant simply truncating the interpolant. From a numerical
point of view this procedure can be accomplished without thinning the data-sites set $X\subset\Omega$, but simply checking 
if the singular values decay under a certain tolerance. This corresponds to solve the linear sistem related to the kernel 
matrix with a (weighted) total least-squares algorithm.
\\\\The dependence of the basis on a singular value decomposition does not allow to produce an adaptive 
algorithm, but forces to compute a full factorization of the matrix for each fixed points distribution. In this sense, 
it would be interesting to adapt our method to work for approximation based on compactly 
supported kernels. Indeed, altough it is possible to use them as any other kernel as done in Section \eqref{sec:newton}, a
 more specific implementation could benefit from the compact support and hence produce sparse kernel matrices. In this 
setting there are eigenvalue algorithms optimized for finding only a small subset of the full spectrum of a matrix, and then
it would be possible to compute an approximant based only on eigenvalues upon a certain tolerance. 

\end{document}